\documentclass[11pt,openany]{article} 
\usepackage{amssymb} \usepackage{amsmath,amscd,mathrsfs}

\usepackage{graphpap}
\usepackage{pstricks}
\usepackage{pst-all}
\usepackage{pstcol}
\usepackage{color}

\usepackage{enumitem}
\usepackage{hyperref}
\usepackage[active]{srcltx}

\newtheorem{thm}{Theorem}
\newtheorem{defin}{Definition}
 \newtheorem{cor}{Corollary}

\newtheorem{prop}{Proposition}

\def\C{{ \! \rm \ I\!\!\!C}}

\def\N{{ \! \rm \ I\!N}}

\def\R{{ \! \rm \ I\!R}}

  \def\det{\rm det}    
   \def \square{\hbox {$\sqcup
    $\llap {$\sqcap $}}} 
  
\newcommand{\parcial}[2]{\frac{\partial#1}{\partial#2}}
\newcommand{\norm}[1]{\vert \vert #1 \vert \vert}
\newcommand{\normv}[1]{ \vert #1 \vert }

\title{Analytic dilation for Laplacians on manifolds with corners of codimension 2} \author{Leonardo A. Cano Garc\'{i}a}
\begin{document}
\maketitle
\begin{abstract}
The analytic dilation method was originally used  in the context
of many body Schr\"odinger operators. In this paper we adapt it to
the context of compatible Laplacians on complete manifolds with
corners of codimension two. As in the original setting of
application we show that the method allows us to: First,   meromorphically extend the matrix elements associated to
analytic vectors. Second, to prove absence of singular spectrum. Third, to find a discrete set that contains the accumulation points
of the pure point spectrum,  and finally, it provides  a
theory of quantum resonances. Apart from these
results, we win also a deeper understanding of the essential
spectrum of compatible Laplacians on complete manifolds with
corners of codimension 2.
\end{abstract}
\section{Introduction}
In the spectral analysis of self-adjoint operators a fundamental problem is to show whether or not the singular spectrum exists and identify the set of accumulation points of the pure point spectrum. In this paper we solve this problem for compatible Laplacians on complete manifolds with
corners of codimension two by adapting the analytic dilation method to this setting.
\\
\\
The method of analytic dilation was originally applied to
$N$-particle Schr\"odinger operators, a classic reference in that
setting is \cite{GERARD}. It has also been  applied to the
black-box perturbations of the Euclidean Laplacian in the series
of papers \cite{SjostrandZworski1} \cite{SjostrandZworski2}
\cite{SjostrandZworski3} \cite{SjostrandZworski4}. In the paper
\cite{B} it is used to study Laplacians on hyperbolic manifolds.
The analytic dilation has also been applied to the study of the
spectral and scattering theory of quantum wave guides and
Dirichlet boundary domains, see e.g. \cite{DUEXMES} \cite{KOVSAC}.
It has also been applied to arbitrary symmetric spaces of
noncompact types in the papers \cite{MV1} \cite{MV2} \cite{MV3}.
In each of these settings new ideas and new methods carry out. In
this paper we develop the analytic dilation method for Laplacians
on complete manifolds with corners of codimension 2.
\\
\\
Let us begin recalling the geometric setting in which our results will be stated. Let $X_0$ be a compact  manifold with boundary $M$. We say that
$X_0$ has a corner of codimension 2 if:
\begin{itemize}
\item[i)] There exists a hypersurface $Y$  of $M$ which divides
$M$ in two manifolds with boundary $M_1$ and $M_2$, i.e. $M=M_1
\cup M_2$ and $Y=M_1 \cap M_2$. \item[ii)] $X_0$ is endowed with a
Riemannian metric $g$ that is a product metric on small
neighborhoods of  the $M_i$'s and the corner $Y$.
\end{itemize}
\begin{center}
\psset{unit=0.5cm}
    \begin{pspicture}(4,1)(12,7)
            \pscurve[](6.5,4)(7,3.7)(7.5,3.5)(8,4)
                \pscurve[](7,3.7)(7.5,3.9)(7.5,3.5)
            \pscurve[](10,6)(7,5.8)(5,5.7)(7,5.5)(9,5.4)
            \pscurve[](10,6)(9.5,3)(9,5.4)
            \pscurve[](5,5.7)(6,2.8)(8.5,3)(9.5,3)
            \rput(7.2,6.5){$M_1$}
            \rput(10.8,4.6){$M_2$}
            \pscircle*(10,6){0.2}
            \pscircle*(9,5.4){0.2}
            \psline(9,5.4)(11.5,6.1)
            \psline(10,6)(11.5,6.1)
            \rput(11.9,6.3){$Y$}
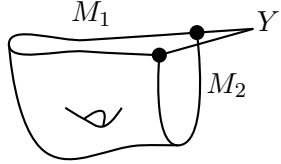
            \rput(8,1.8){Figure 1. Compact manifold with corner of codimension 2}
        \end{pspicture}
\end{center}
We construct from $X_0$ a complete manifold $X$ by attaching
$([0,\infty) \times M_1 )$,  $([0,\infty) \times M_2 )$ and
filling the rest with $( [0,\infty) \times  [0,\infty) \times Y
)$. As a set,
\begin{equation}
 X:=X_0 \cup ([0,\infty) \times M_1 ) \cup ([0,\infty) \times M_2 ) \cup ( [0,\infty) \times [0,\infty) \times Y ),
\end{equation}
and it has the natural differential structure and Riemannian
metric that are compatible with the product structures at the
boundary of $X_0$. The manifold $X$ has associated a natural
exhaustion given by:
\begin{equation}\label{eq:def exhaust}
 X_T:=X_0 \cup ([0,T] \times M_1 ) \cup ([0,T] \times M_2 ) \cup ( [0,T]|^2 \times Y
 ).
\end{equation}
\vspace{0.3cm}
\begin{center} \psset{unit=0.3cm}
    \begin{pspicture}(0,0)(24,10)
            \pscurve[](6.5,4)(7,3.7)(7.5,3.5)(8,4)
            \pscurve[](7,3.7)(7.5,3.9)(7.5,3.5)
        \pscurve[](5,5.7)(6,2.8)(8.5,3)(9.5,3)
        \pscurve[](10,11)(7,10.8)(5,10.7)(7,10.5)(9,10.4)
        \psline(9,10.4)(14,10.4)
        \psline(14,10.4)(14,5.3)
        \psline(9.5,3)(14.5,3)
        \psline(5,5.7)(5,10.7)
        \psline[linestyle=dashed](10,11)(15,11)
        \psline[linestyle=dashed](15,11)(15,6)
        \pscurve[](15,6)(14.5,3)(14,5.4)
        \rput(9,1){Figure 2. $X_T$, element of the exhaustion of $X$.}
        \psline(9,5.4)(9,10.4)
        \psline(5,5.7)(5,10.7)
        \psline[linestyle=dashed](10,6)(10,11)
        \rput(6.8,8.5){{\tiny $[0,T] \times M_1$}}
        \pscurve[](10,6)(7,5.8)(5,5.7)(7,5.5)(9,5.4)
        \pscurve[](10,6)(9.5,3)(9,5.4)
        \pscurve[](5,5.7)(6,2.8)(8.5,3)(9.5,3)
        \pscurve[](15,6)(14.5,3)(14,5.4)
        \psline(9.5,3)(14.5,3)
        \psline[linestyle=dashed](15,6)(10,6)
        \psline(14,5.4)(9,5.4)
        \rput(11.5,4.7){{\tiny$ [0,T] \times M_2$}}
        \rput(12,8.5){{\tiny $ [0,T]^2\times Y$}}
        \rput(6,4.5){{\tiny $X_0$}}
\end{pspicture}
\end{center}
For each $T \in [0,\infty)$, $X$ has two submanifolds with
cylindrical ends, namely $M_i \times \{T\} \cup (Y \times \{T\})
\times [0,\infty) $, for $i=1,2$. We denote these manifolds by
$Z_i$.
\begin{center}
      \psset{unit=0.5cm}
          \begin{pspicture}(0,-2)(10,10)
                  \psline(0,0)(0,10)
                   \qline(0,0)(10,0)
        \psline(2,0)(2,10)
        \psline(0,2)(10,2)
        \rput(1,1){$X_T$}
        \rput(9,9){$[T,\infty)^2 \times Y$}
        \rput(9,2.8){$Z_1 \times \{T\}$}
        \rput(2.6,10.5){$Z_2\times \{T\}$}
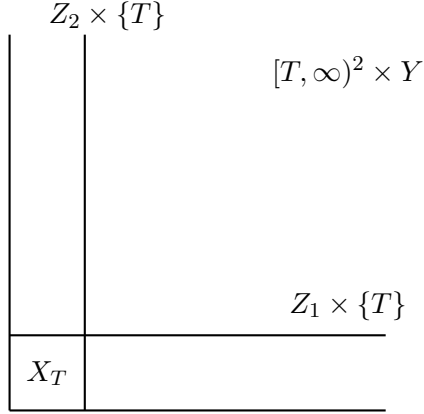
        \rput(7,-1){Figure 3. Sketch of a complete manifold with corner of codimension 2}
    \end{pspicture}
\end{center}
Let us now consider the operator $\Delta:=d^*d:C^\infty_c(X)\to
L^2(X)$, in local coordinates:
\begin{equation}
\Delta=\frac{-1}{\sqrt{\normv{\det (g)}}}\parcial{}{x_i}\left(g^{ij}\sqrt{\normv{\det (g)}}\parcial{}{x_j} \right).
\end{equation}
The operator $\Delta$ is essentially self-adjoint and we denote by
$H$ its self-adjoint extension. In section \ref{subsec:Manifolds
with corners} we shall consider compatible  Laplacians.
\\
\\
For $i=1,2,$ since $Z_i$ is a complete manifold, the Laplacian
$\Delta:C^\infty_c(Z_i) \to L^2(Z_i)$ is essentially self-adjoint;
we denote by $H^{(i)}$ its self-adjoint extension. Similarly
$\Delta:C^\infty(Y) \to L^2(Y)$ is  essentially  self-adjoint, and
we denote its self-adjoint extension by $H^{(3)}$. The analytic
dilation of a many-body Schr\"odinger operator depends on the
analytic dilation of its subsystem Hamiltonians. In a similar way
the analytic dilation of $H$ is described in terms of the spectral
theory of the operators $H^{(1)}$, $H^{(2)}$ and $H^{(3)}$.
\\
\\
For $\theta>0$, the operator $U_{i,\theta}:L^2(Z_i) \to L^2(Z_i)$
is essentially  the dilation operator by $\theta+1$ up to a
compact set. More precisely:
\begin{equation}
U_{i,\theta}f(x)=\begin{cases}
f(x)& \text{ for } x \in M_i.\\
\\
(\theta+1)^{1/2}f((\theta+1)u,y)& \text{ for } x=(u,y) \in [0,\infty) \times Y\\
&\hspace{0.5cm}\text{ and for  } u \text{ big enough},
\end{cases}
\end{equation}
and $U_{i,\theta}f$ is extended to the whole $Z_i$ in such a way that it sends $C^\infty_c(Z_i)$ into $C^\infty_c(Z_i)$, and it becomes a unitary operator on $L^2(Z_i)$. Details will be worked out in section \ref{sec:analytic dil. cyl end}. Similarly, the operators $U_\theta:L^2(X) \to L^2(X)$ are defined by
\begin{equation}
U_{\theta}f(x)=\begin{cases}
f(x)& \text{ for } x \in X_0.\\
\\
(\theta+1)^{1/2}U_{i,\theta}f((\theta+1)u_i, z_i)& \text{ for } x=(u_i,z_i) \in [0,\infty) \times  Z_i\\
&\hspace{0.5cm} \text{ and for  } u_i \text{ big enough}.
\end{cases}
\end{equation}
Again $U_\theta f$ is extended to the whole $X$ in such a way
that, for $f \in C^\infty_c(X)$, $U_\theta f \in C^\infty_c(X)$,
and  $U_\theta$ becomes  a unitary operator in $L^2(X)$. Details
will be given in section \ref{def:Utheta corner}.
\\
\\
For $\theta \in [0,\infty)$, define $H_\theta:=U_\theta H
U_\theta^{-1}$, a closed operator with domain
\begin{equation}
\mathscr{W}_{2}(X):=\{f \in L^2(X): \Delta_{dist}f \in L^2(X)\},
\end{equation}
the second Sobolev space associated to $(X,g)$. We define the set:
\begin{equation}\label{def: Gamma}
\begin{split}
\Gamma:=\{\theta:=\theta_0+i\theta_1 \in \C:\theta_0>0, \theta_0
\geq \normv{\theta_1} \text{ and } Im(\theta)^2<1/2\}.
\end{split}
\end{equation}
We will extend the family $H_\theta$ from $[0,\infty)$ to
$\Gamma$.
\begin{center}
\psset{unit=0.5cm}
    \begin{pspicture}(-3,-2)(12,11)
                \pspolygon[fillstyle=hlines](12,0)(8,0)(4,4)(8,8)(12,8)
            \psline{->}(-1,4)(12,4)
            \psline{->}(4,-1)(4,10)
            \rput(3.2,8){$\frac{\sqrt{2}}{2}$}
            \rput(3,0){$-\frac{\sqrt{2}}{2}$}
            \rput(8,3.5){$\frac{\sqrt{2}}{2}$}
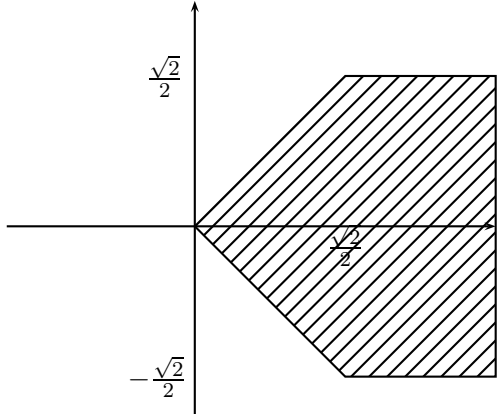
            \rput(4,-1.5){Figure 4. Sketch of the region $\Gamma$}
        \end{pspicture}
\end{center}
We prove:
\begin{thm}\label{thm:Afamily} The family $(H_\theta)_{\theta \in [0,\infty)}$ extends to an holomorphic
family for $\theta \in \Gamma$,  which satisfies:
\\
\\
1) $H_\theta$ is a closed operator with domain
$\mathscr{W}_{2}(X)$ $\theta \in \Gamma$.
\\
\\
2) For $\varphi \in \mathscr{W}_{2}(X)$ the map $\theta \mapsto
H_\theta \varphi$ is holomorphic in $\Gamma$.
\end{thm}
An holomorphic family of operators  satisfying (1) and (2) will be
called {\bf a holomorphic family of type A}. This theorem is
proved using the analogous result that the family
$\{H^{(i),\theta}\}_{\theta \in [0,\infty)}$    extends to a
holomorphic family of type A in $\Gamma$, where $H^{(i),\theta}$
denotes the closed operator associated to $U_{i,\theta}
\Delta_{Z_i} U_{i,\theta}^{-1}$ with domain
\begin{equation}
\mathscr{W}_{2}(Z_i):=\{f \in L^2(Z_i):\Delta_{dist}(f) \in L^2(Z_i)\},
\end{equation}
the second Sobolev space associated to $(Z_i,g_i)$.
\\
\\
The families $H_\theta$ and $H^{(i)}_\theta$ extend to sets larger
than $\Gamma$, but $\Gamma$ is enough for our outlined goals.  In
particular, we choose the domain $\Gamma$ because for $\theta \in
\Gamma$ we can prove that $H^{(i)}_\theta$ is $m$-sectorial (see
section \ref{sec: Deltatheta is sectorial}).
We define
\begin{equation}
\theta':=\frac{1}{ (\theta+1)^2}.
\end{equation}
The parameter $\theta'$ is very important in the description of
the essential spectrum of $H_\theta$ as we can see in  the next
theorem that will be proved in section \ref{section:calculating
essent. spec}.
\begin{thm}\label{thm: ess Htheta}
For $\theta \in \Gamma$,
\begin{equation}
\begin{split}
\sigma_{ess}(H_\theta)=&\bigcup_{\mu \in \sigma(H^{(3)})} (\mu+\theta'[0,\infty))\\
&\cup  \bigcup_{\lambda_1 \in \sigma_{pp}(H^{(1),\theta})} \left(\lambda_1 +\theta'[0,\infty) \right)\\
&\cup  \bigcup_{\lambda_2 \in \sigma_{pp}(H^{(2,\theta)})} \left(\lambda_2 +\theta'[0,\infty) \right).
\end{split}
\end{equation}
\end{thm}
In section \ref{sect:analytic vect corners},  we associate to
$(U_\theta)_{\theta \in [0,\infty)}$ a set $\mathscr{V} \subset
\mathscr{W}_{2}(X)$   that satisfies:
\begin{itemize}
\item[i)]$\mathscr{V}$ is dense in $L^2(X)$. \item[ii)] for
$\varphi \in \mathscr{V}$, $U_\theta \varphi $ is defined for all
$\theta \in \Gamma$. \item[iii)] $U_\theta \mathscr{V}$ is dense
in $L^2(X)$ for all $\theta \in \Gamma$.
\end{itemize}
The elements of a subset of $\mathscr{W}_{2}(X)$  which satisfies
i) and ii) will be called {\bf analytic vectors}.  We denote by
$\Lambda$ the left-hand plane, more explicitly:
\begin{equation}\label{eq: def Lambda}
\Lambda:=\{(x,y) \in \C: x < 0\}.
\end{equation}
We denote by $R(\lambda)$ the resolvent of $H$ and by $R(\lambda)$
the resolvent of $H_\theta$.  Using the general analytic dilation
theory of Aguilar-Balslev-Combes (see \cite{B}) we describe the
nature of the spectrum of $H$. This theory is based on:
\begin{itemize}
\item[i)]  The knowledge of the essential spectrum of $H_\theta$,
provided by theorem  \ref{thm: ess Htheta}. \item[ii)] The
following equation, that is consequence of the unitarity of
$U_\theta$,
\begin{equation}\label{eq:ext resol elem}
\left<R(\lambda)f,g
\right>_{L^2(X)}=\left<R(\lambda,\theta)U_\theta f,U_\theta g
\right>_{L^2(X)},
\end{equation}
for $f,g \in \mathscr{V}$ and $\theta \in [0,\infty)$.
\end{itemize}
Since the right-hand side of (\ref{eq:ext resol elem}) is defined
for $\lambda \in \Lambda$ and $\theta \in \Gamma$, (\ref{eq:ext
resol elem}) provides a meromorphic extension of $\lambda \mapsto
\left<R(\lambda)f,g \right>_{L^2(X)}$ from $\Lambda$ to
$\C-\sigma(H_\theta)$. From this, we deduce the following theorem.
\begin{thm}
1) For $f,g \in \mathscr{A}$ the function $\lambda \mapsto \langle
R(\lambda)f,g \rangle_{L^2(X)}$ extends from $\Lambda$ to
$\C-\sigma(H_\theta)$.
\\
\\
2) For all $\theta \in \Gamma$, $H_\theta$ has no singular
spectrum.
\\
\\
3) The accumulation points of $\sigma_{pp}(H)$ are contained in
$\{\infty\}\cup \sigma(H^{(3)}) \cup  \cup_{i=1}^2
\sigma_{pp}(H^{(i)})$.
\end{thm}
The meromorphic extension of the resolvent entries can be used to
extend generalized eigenfunctions that describe natural wave
operators whose image is the complete set of absolute continuous
states associated to $H$. These results are worked out in detail
in \cite{CANOGENEIGEN} \cite{CANOASYMPCOMPLETENESS}
\cite{CANOTHESIS}.
\\
\\
In appendix \ref{chap:DHSV}, based on \cite{DHSV}, we
geometrically refine the notion of singular  spectrum introducing
 what we call boundary Weyl sequences; we use such
refinament to compute the essential spectrum in section
\ref{section:calculating essent. spec}. In appendix
\ref{chap:Ichinose} we give a brief introduction to sectorial
operators in such a way that we can enunciate the Ichinose lemma
(see theorem \ref{thm:Ichinoselemma} in appendix
\ref{chap:Ichinose}). Given two closed operators $A$ and $B$
acting on Hilbert spaces $\mathscr{H}_1$ and $\mathscr{H}_2$,  the
Ichinose lemma  gives sufficient conditions for having
$\sigma(A\otimes 1 +1\otimes B)=\sigma(A)+\sigma(B)$. It is
important for our results because, intuitively, the operator
$H_\theta$ is of the form $A\otimes 1 +1\otimes B$ at infinity.
\\
\\
This paper is based on the PhD thesis of the author, realized
under the supervision of Werner M\"uller, at the University of Bonn,
the author is grateful with him for his support. Thanks are also due to  Rafe Mazzeo for helping to contextualize this document, and
Julie Rowlett and Alexander Cardona for their help with its final
presentation.
\section{Analytic dilation on complete manifolds with cylindrical end}
\label{sec:analytic dil. cyl end}In this section we generalize the
method of analytic dilation to complete  manifolds with
cylindrical end. The results of this section are consequence of
\cite{KALVIN1} and they are expected from the results of analytic
dilation in wave guides (see \cite{DUEXMES} \cite{KOVSAC}).
Because of that, and since the approach in the proof of the main
results of this section is similar to the approach in section
\ref{subsec:Manifolds with corners}, we give most of the results
without proof. The interested reader is refered to
\cite{CANOTHESIS} or \cite{KALVIN1} for further details.
\\
\\
The analytic dilation on complete manifolds with cylindrical  end
will be important in section \ref{subsec:Manifolds with corners}
because the analytic dilation of $H$ is described in terms of the
analytic dilations of $H^{(1)}$ and  $H^{(2)}$, compatible
Laplacians on $Z_1$ and $Z_2$. In fact, one of our main results,
theorem \ref{thm: ess Htheta}, shows that  the essential spectrum
of $H_\theta$ is described in terms of the pure point spectrum of
$H^{(1)}_\theta$ and $H^{(2)}_\theta$, the dilated Laplacians
associated to $H^{(1)}$ and  $H^{(2)}$.
\subsection{Manifolds with cylindrical end  and their compatible Laplacians}
\label{sect:analyticdilationcyl}Let  $Z_0$ be a compact Riemannian
manifold with boundary $Y:=\partial Z_0$. We say that $Z_0$ is a
{\bf compact manifold with cylindrical end} if there exists a
neighborhood, $Y \times (-\epsilon,0]$, of the boundary $Y$ such
that Riemannian metric of $Z_0$ is a product metric i. e. a metric
of the form $g_Y+du\otimes du$ where $g_Y$ is a Riemannian metric
on $Y$ and $u$ is the variable on $(-\epsilon,0]$.
\begin{center}
\psset{unit=0.5cm}
    \begin{pspicture}(0,0)(24,7)
         \rput(0.8,7){$Z_0$} \pscurve[](2,2)(0,4)(2,6)
         \pscurve[](2,6)(3,6)(4,5.9)(6,5.8)
         \psline(6,5.8)(8,5.8) 
         \pscurve[](2,2)(3,2)(4,2.1)(6,2.2)
         \psline(6,2.2)(8,2.2)
         \pscurve[](1.5,4)(2,3.7)(2.5,3.5)(3,4)
        \pscurve[](2,3.7)(2.5,3.9)(2.5,3.5) \psellipse[](8,4)(0.7,1.8)
        \psellipse[linestyle=dashed](6.8,4)(0.7,1.8)
        \psellipse[](15,4)(0.7,1.8)
        \rput[r](14.1,4){$Y:=\partial Z_0=$}
        \psline{|->}(6.2,4)(7,4)
        \rput(6.75,3.5){$\partial u$}
        \rput(10,1){Figure 5. Compact
manifold with cylindrical end.}
\end{pspicture}
\end{center}
We make from $Z_0$ a  complete manifold $Z$ by attaching the
infinite cylinder $Y \times [0,\infty)$ to $Z_0$. We have then:
\begin{equation}\label{eq:man cyl end}
Z:=Z_0 \cup_Y (Y \times [0,\infty)),
\end{equation}
where we are identifying the boundary of $Z_0$ with $Y \times
\{0\}$. We extend the smooth structure and the Riemannian metric
naturally. The manifold $Z$ is called {\bf complete manifold with
cylindrical end}. It looks as follows:
\begin{center}
\psset{unit=0.5cm}
\begin{pspicture}(0,0)(24,7)
\pscurve[](2,2)(0,4)(2,6) \pscurve[](2,6)(3,6)(4,5.9)(6,5.8)
\psline(6,5.8)(20,5.8) \pscurve[](2,2)(3,2)(4,2.1)(6,2.2)
\psline(6,2.2)(20,2.2) \pscurve[](1.5,4)(2,3.7)(2.5,3.5)(3,4)
\pscurve[](2,3.7)(2.5,3.9)(2.5,3.5) \psellipse[](20,4)(0.7,1.8)
 \rput(10,1){Figure 6. Complete manifold with cylindrical
end.}
\end{pspicture}
\end{center}
Let $E$ be a vector bundle over $Z$ with an Hermitian metric. We
assume  that there exists $E'$ an Hermitian vector bundle over $Y$
such that $E \vert_{Y \times [0,\infty)}$ is the pull back of $E'$
by the projection $\pi:Y \times \R_+ \to Y$. We suppose that the
Hermitian metric of $E$ is the pullback of the Hermitian metric of
$E'$.  Let $\Delta$ be a  generalized Laplacian on $Z$, i.e.
$\sigma_2(\Delta)(z,\xi)= \normv{\xi}_{g_z}^2$. We assume
furthermore that on $Y \times [0,\infty)$
\begin{equation} \label{eq: generalized laplacian}
\Delta=-\parcial{^2}{u^2}+\Delta_Y,
\end{equation}
where $\Delta_Y$ is a generalized Laplacian acting on $C^\infty(Y,E')$. In fact, we will denote by $\Delta_Y$ the operator acting on distributions and the self-adjoint operator induced by $(\Delta_Y, C^\infty(Y,E'))$.
\\
\\
A  Laplacian satisfying the previous assumptions is called a {\bf
compatible  Laplacian}.
\subsection{The definition of $U_\theta$}
\label{analytic dilation cyl end}Let $0<K<R$ and $\varphi \in C^\infty(\R)$ with $\varphi'\geq 0$ such that:
\begin{equation}
\varphi(u):=\begin{cases}
0& \text{for } 0<u<K\\
1 & \text{for } R<u<\infty.
\end{cases}
\end{equation}
Let $\theta \in [0,\infty)$, define the function:
\begin{equation} \label{eq:calpsitheta1}
\psi_\theta(u):=(\varphi(u)\theta+1)u=\varphi(u)u\theta+u,
\end{equation}
for $u\in \R_+$. Observe that
\begin{equation*}
\psi_\theta(u)=\begin{cases}
u & u<K\\
(\theta+1)u &\text{ for }u>R.
\end{cases}
\end{equation*}
We calculate the first derivatives of $\psi_\theta$ with respect
to $u$:
\begin{equation} \label{eq:calpsitheta2}
\psi_\theta'(u):=\parcial{}{u}(\psi_\theta)(u)=\varphi'(u)u\theta+\varphi(u)\theta+1.
\end{equation}
\begin{equation}\label{eq:calpsitheta3}
\psi_\theta''(u):=
\parcial{^2}{u^2}\psi_\theta(u)=
\varphi''(u)u\theta+2\varphi '(u)\theta.
\end{equation}
and
\begin{equation} \label{eq:calpsitheta4}
\psi_\theta'''(u):=
\parcial{^3}{u^3}\psi_\theta(u)=
\varphi'''(u)u\theta+3\varphi''(u)\theta.
\end{equation}
We define $U_\theta:L^2(Z,E) \to L^2(Z,E)$:
\begin{equation}
U_\theta f(z)=\begin{cases}f(z)& \text{ for } z \in Z_0\\
f(y, \psi_\theta(u))\psi_\theta'(u)^{1/2}&\text{ for } z=(y,u) \in Y \times \R_+.
\end{cases}
\end{equation}
Observe that, for $\theta>0$, the function $\psi_\theta$ is invertible (because $\psi_\theta'(u)\geq 1$ for $u\geq 0$). We will denote its inverse by $\alpha_\theta$.
\\
\\
For $\theta \in \R_+$ a natural inverse of $U_\theta$ is given by:
\begin{equation}
U_\theta^{-1} f(z):=\begin{cases}f(z_0), \text{ for } z=z_0 \in Z_0.\\
f(y,\alpha_\theta(u))\psi_\theta'(\alpha_\theta(u))^{-1/2} \text{ for } (y,u) \in Y \times \R_+.
\end{cases}
\end{equation}
We observe that for $f \in C^\infty(Z,E)$, $U_\theta f$ belongs to
$ C^\infty(Z,E)$ and, if $f \in C^\infty_c(Z,E)$, then $U_\theta f
\in C^\infty_c(Z,E)$. It is easy to see:
\begin{prop}
For $\theta \in (0,\infty)$, $U_\theta$ induces a unitary operator acting on $L^2(Z,E)$.
\end{prop}
In the next section we extend the family of operators $(U_\theta
\Delta U_{\theta}^{-1})_{\theta \in (0,\infty)}$ to parameters
$\theta \in \Gamma$ (see (\ref{def: Gamma})).
\subsection{The family $\Delta_\theta$}
\label{sec:family of type A cylend}Calculating explicitly $U_\theta \parcial{^2}{u^2}U_\theta^{-1}f(y,u)$, for $(y,u) \times Y \times  \R_+$, we obtain:
\begin{equation}\label{eq.calculation of Deltatheta}
\begin{split}
&\Delta_\theta f(y,u)= \Delta_Y f(y,u)-\parcial{^2}{u^2}f(y,u)(\alpha_\theta '(\psi_\theta(u)))^2\\
&\hspace{0.5cm}-\parcial{}{u}f(y,u)\alpha_\theta ''(\psi_\theta(u))+\parcial{}{u}f(y,u) (\psi_\theta'(u))^{-1} \psi_\theta''(u) (\alpha_\theta'(\psi_\theta(u)))^2\\
&\hspace{0.5cm}-3/4 f(y,u)\psi_\theta'(u)^{-2} (\psi_\theta''(u))^2 (\alpha_\theta'(\psi_\theta(u)))^2+1/2 f(y,u) (\psi_\theta'(u))^{-1} \psi_\theta'''(u) (\alpha'(\psi_\theta(u)))^2\\
&\hspace{0.5cm}+1/2 f(y,u) (\psi_\theta'(u))^{-1} \psi_\theta''(u)
\alpha_\theta''(\psi_\theta(u)).
\end{split}
\end{equation}
Observe that, on $Y \times \R_+$, we have:
\begin{equation} \label{eq:coeffiecientdeltatheta}
\Delta_\theta=a_2(\theta,u)\parcial{^2}{u^2}+a_1(\theta,u)\parcial{}{u}+a_0(\theta,u)+\Delta_Y,
\end{equation}
where $a_2(\theta,u), a_1(\theta,u)$ and $a_0(\theta,u)$ are given
by
\begin{equation}\label{eq:a1 a2 a3 theta}
\begin{split}
&a_2(\theta,u):=(\alpha'_\theta(\psi_\theta(u)))^2=\frac{-1}{(\psi_\theta'(u))^2};\\
&a_1(\theta,u):=-\alpha_\theta''(\psi_\theta(u))+(\psi_\theta'(u))^{-1}\psi_\theta''(u)\alpha_\theta'(\psi_\theta(u))^2;\\
&a_0(\theta,u):=1/2\psi_\theta'(u)^{-1}\psi_\theta'''(u)\alpha_\theta'(\psi_\theta(u))^2\\
&\hspace{1.4cm}+1/2\psi_\theta'(u)^{-1}\psi_\theta''(u)\alpha_\theta''(\psi_\theta(u))\\
&\hspace{1.4cm}-3/4\psi_\theta'(u)^{-3/2}\psi_\theta''(u)^2\alpha_\theta'(\psi_\theta(u))^2.
\end{split}
\end{equation}
Notice that, for all $u \in  \R_+$, $a_2(\theta,u), a_1(\theta,u)$
and $a_0(\theta,u)$ are well defined and holomorphic for
$Re(\theta)\geq 0$. We remark also that  $$a_k(\theta,.) \in
C^\infty(\R_+)$$ for $k=0,1,2$ and $Re(\theta)>0$.  We will
continue  denoting  $a_2(\theta,u), a_1(\theta,u)$ and
$a_0(\theta,u)$, the coefficients of $\parcial{^2}{u^2}$,
$\parcial{}{u}$ and $Id$, respectively,  for the operator
$\Delta_\theta$ localized in $Y\times \R_+$.
\\
\\
From now on, given $\theta \in \C-(-\infty,0)$, we define
$\theta'$ by
\begin{equation}
\theta':= \frac{1}{(\theta+1)^2}.
\end{equation}
The parameter $\theta'$ will appear naturally in the description
of $\sigma_{ess}(\Delta_\theta)$ (see equation
(\ref{eq:essspec1})). The next proposition follows easily  from
(\ref{eq.calculation of Deltatheta}).
\begin{prop}\label{prop: formula htheta cylend}
Let $f \in C^\infty(Z,E)$. For $Re(\theta) \geq 0$, the formula for $\Delta_\theta$ reduces for $(y,u) \in  Y \times (0,K)$ to:
\begin{equation}
\Delta_\theta f(u,y)=-\parcial{^2}{u^2}f(u,y)+\Delta_Y f(u,y);
\end{equation}
and, for $(y,u) \in  Y \times (R,\infty)$ to:
\begin{equation}
\Delta_\theta f(u,y)=-\theta' \parcial{^2}{u^2}f(u,y)+\Delta_Y f(u,y).
\end{equation}
\end{prop}
The next proposition is a technical tool that can be deduced from
(\ref{eq:a1 a2 a3 theta}).
\begin{prop}\label{prop:aibounded}
Let $a_0(\theta,u)$, $a_1(\theta,u)$ and $a_2(\theta,u)$ be given
by (\ref{eq:a1 a2 a3 theta}). If $\vert \theta \vert < N$ and
$Re(\theta)\geq 0$, then there exists a $C(N) \in \R_+$
independent of $\theta$ and $u\in \R_+$ such that for $i=0,1,2$:
\begin{equation}
\vert a_i(\theta,u)\vert \leq C(N),
\end{equation}
and,
\begin{equation}
\vert \parcial{}{\theta_i}\left(a_i\right)(\theta,u)\vert \leq C(N),
\end{equation}
where $\theta:=\theta_1+i\theta_2$.
\end{prop}
We recall some definitions and results on manifolds with bounded
geometry and their natural vector bundles, references for them are
\cite{EICHHORN} and \cite{Shubin}. A Riemannian manifold { \bf $M$
has bounded geometry} if its injectivity radious is positive and
if all the derivatives of the curvature tensor, in the geodesic
coordinates, are uniformly bounded. Let $D$ be a linear
differential operator in $Diff(E,F)$, $E$ and $F$ vector bundles
over $M$ with bounded covariant derivatives $\nabla^E$ and
$\nabla^F$ over $M$, a manifold with bounded geometry. {\bf $D$
has bounded coefficients} if, in geodesic coordinates and in any
synchronous maps, all the derivatives are uniformly bounded. It is
easy to see that the manifold $Z$ given by (\ref{eq:man cyl end})
has bounded geometry and the compatible Laplacians in $Diff(E)$
have bounded coefficients. A differential operator $D$ of order
$m$ {\bf is a uniformly elliptic operator}, if the following
inequality holds:
\begin{equation}
\normv{\sigma_{m}(D)^{-1}(z,\xi)}^{-1} \leq \normv{\xi}^{-m},
\end{equation}
for $z \in Z$ and $\xi \in T_zZ$. If $A$ is a second order
differential uniformly elliptic operator, then the norm  $f
\mapsto \norm{f}_{L^2(Z,E)} +\norm{A f}_{L^2(Z,E)} $ on
$C^\infty_c(Z,E)$ is equivalent to the norm $f \mapsto
\norm{f}+\norm{\Delta (f)}$. In general, if $A$ is a uniformly
elliptic differential operator of order $m$ acting on a manifold
with bounded geometry, then the norms $f \mapsto \norm{f}+\norm{A
(f)}$ and the $m$-Sobolev norm $\norm{\cdot}_m$, defined naturally
using geodesic coordinates and synchronous frames, are equivalent.
Since $Z$ is a complete manifold and $\Delta:C^\infty_c(Z,E) \to
L^2(Z,E)$ is uniformly elliptic, then $\Delta$ is essentially
self-adjoint (see \cite{Shubin}). Abusing of the notation we
denote by $\Delta$ the differential operator acting on
distributions and the self-adjoint operator itself.
 Denote by
$\mathscr{W}_2(Z,E)$  the closure of $C_c^\infty(Z,E)$ with
respect to the norm $\norm{f}_2:=\norm{f}+\norm{\Delta f}$ for $f
\in C_c^\infty(Z,E)$. We call $\mathscr{W}_2(Z,E)$ the {\bf second
Sobolev space}.
\\
\\
Using the theory of manifolds with bounded geometry sketched above
it is straight to prove the next theorem.
\begin{thm}\label{thm:analyextendelta}
The family $(\Delta_\theta)_{\theta \in \R_+}$ extends to an analytic family of type A for $Re(\theta) > 0$  i.e.
\\
\\
i) $\Delta_\theta$ are closed operators with $Dom(\Delta_\theta)$ independent of $\theta$. More precisely, $Dom(\Delta_\theta)=\mathscr{W}_2(Z,E)$.
\\
\\
ii) For every $f \in \mathscr{W}_2(Z,E)$ the map $\theta \mapsto \Delta_\theta f$ is analytic for $Re(\theta) > 0$.
\end{thm}
\subsection{The essential spectrum of $\Delta_\theta$}
\label{section:essential spectrum cyl end} Recall that, given a
closed operator $A$, the pure point spectrum, discrete spectrum,
and essential spectrum are the sets given by
\begin{equation}\label{eq: def ess d pp}
\begin{split}
&\sigma_{pp}(A):=\{\lambda \in \C: \text{ is an eigenvalue of $A$}\},\\
&\sigma_{d}(A):=\{\lambda \in \C: \text{ $\lambda$ is an isolated eigenvalue of $A$ of finite algebraic multiplicity}\},\\
&\sigma_{ess}(A):=\sigma(A)-\sigma_{d}(A),
\end{split}
\end{equation}
respectively. Our next goal is to prove the equality:
\begin{equation} \label{eq:essspec1}
\sigma_{ess}(\Delta_\theta)=\bigcup_{i=0}^\infty \left(\mu_i+\theta' [0,\infty)\right),
\end{equation}
where $\sigma(\Delta_Y):=\{\mu_i\}_{i=0}^\infty$. The first step
towards (\ref{eq:essspec1}) is to prove:
\begin{equation} \label{eq:Ness equal F theta}
N_{ess}(\Delta_\theta)=\bigcup_{i=0}^\infty \left(\mu_i+\theta' [0,\infty)\right),
\end{equation}
where $N_{ess}$ is the set defined in appendix  \ref{chap:DHSV},
definition \ref{def: Ness}. Theorem \ref{thm:DHSVNess} and
equation (\ref{eq:Ness equal F theta}) imply
\begin{equation}
\sigma_{ess}(\Delta_\theta)=N_{ess}(\Delta_\theta),
\end{equation}
and hence (\ref{eq:essspec1}).
\\
\\
The proof of equation (\ref{eq:Ness equal F theta}) is based on
the manipulation of singular sequences (see definition
\ref{def:singseq} in appendix \ref{chap:DHSV}). In
\cite{CANOTHESIS} we prove that singular sequences associated to
$\Delta_\theta^i$, the closed operator associated to $-\theta'
\parcial{^2}{u^2}+\mu_i$ with Dirichlet boundary conditions, induce singular sequences associated
$\Delta_\theta$; and the other way around, singular sequences
associated to $\Delta_\theta$ induce singular sequences associated
$\Delta_\theta^i$ for some $i$. A fundamental tool for formalizing
these ideas is the Rellich theorem. As we said in the introduction
we do not give details here since a similar approach will be used
in section \ref{section:calculating essent. spec}.
\begin{thm}\label{thm:gluecyl}
1) Let $Re(\theta)\geq 0$ and let $(g_n)_{n \in \N}$ be an orthonormal singular sequence associated to $\lambda$ and $-\theta' \parcial{^2}{u^2}+\mu_i$. Then, there exists a subsequence of  $h_n:=(\kappa g_n \phi_i).\norm{\kappa g_n \phi_i}^{-1}$ that is a singular sequence  associated to the operator $\Delta_{\theta}.$
\\
\\
2) Let $Re(\theta)\geq 0$ and let $g_n$ be an orthonormal singular
sequence associated to the operator $\Delta_\theta$ and the value
$\lambda$. Then, there exists $i \in \N$ and a subsequence $s$ of
$\N$ such that the function $u \mapsto \left<\kappa
g_{s(n)}(u,\cdot), \phi_i \right>_{L^2(Y)}$, in $C^\infty(\R_+)$,
is a singular sequence of $\Delta_\theta^i$  and the value
$\lambda$.
\end{thm}
From theorem \ref{thm:gluecyl} follows that:
\begin{equation} \label{eq: Ness equal Ftheta cyl end}
N_{ess}(\Delta_\theta)=\bigcup_{i=0}^\infty
\left(\mu_i+\theta'[0,\infty)\right),
\end{equation}
which, together with theorem \ref{thm:DHSVNess}, imply that
\begin{equation}
\sigma_{ess}(\Delta_\theta)=N_{ess}(\Delta_\theta).
\end{equation}
\subsection{The analytic vectors of $U_\theta$}
\label{analytic vectors cyl end} In this section we construct a
subset, $\mathscr{V} \subset L^2(Z,E)$ such that
\begin{itemize}
\item[i)] $\mathscr{V}$ is a dense subset of $L^2(Z,E)$.
\item[ii)] For $ f \in \mathscr{V}$ the function $\theta \mapsto U_\theta f \in L^2(Z,E)$ makes sense for $\theta \in \C$, $Re(\theta)>0$.
\item[iii)] $U_\theta \mathscr{V}$ is dense in $L^2(Z,E)$ for $Re(\theta)>0$.
\end{itemize}
Recall that we denote by $(\phi_i,\mu_i)_{i=1}^\infty$  a spectral resolution of the operator $\Delta_Y$. Let $\kappa \in C^\infty(\mathbb{R}_+)$ be a function satisfying $0 \leq \kappa \leq 1$, $\kappa'\geq 0$ and
\begin{equation}
\kappa (u)= \begin{cases}
1 & K \leq u< \infty .\\
0 & 1<u \leq K-1.
\end{cases}
\end{equation}
We extend $\kappa$ to $Y\times \R_+$ defining
$\kappa(y,u):=\kappa(u)$ for $(y,u) \in Y \times \R_+$. Making
$\kappa$ equal to $0$ out of its support in $Y \times \R_+$, we
extend $\kappa$ to $Z$. Define the set $\mathscr{P}$ of elements
$h \in L^2(Y \times \R_+,E)$ such that $h$  has a Fourier
expansion of the form
 $h(y,u)=1/u^2\sum_{i=0}^\infty p_i(1/u)\phi_i(u)$, where $p_i(x) \in \mathbb{C}[x]$. Define:
\begin{equation}\label{eq:analytvectors1}
\begin{split}
&\mathscr{V}:=\{(1-\kappa) g+\kappa h: g \in L^2(Z,E) \text{ and } h \in \mathscr{P}\}.
\end{split}
\end{equation}
i) and iii) are consequence of the Stone-Weirstrass theorem.
\subsection{Consequences of Aguilar-Balslev-Combes theory}
In this section we provide the description of
$\sigma_d(\Delta_\theta)$ and $\sigma_{pp}(\Delta_\theta)$ that
the analytic dilation provides. Most of the results that we
compile in the next theorem are consequences of the
Aguilar-Balslev-Combes theory as explained in the book
\cite{HISLOPSIGAL}.
\begin{thm}\label{thm: descriptio delta theta cyl}We have:
\\
\\
a) The set  of non-threshold eigenvalues\footnote{The set of
thresholds of $\Delta$, $\tau(\Delta)$, is by definition
$\sigma(\Delta_Y)$} of $\Delta$ is equal to
$\sigma_{d}(\Delta_\theta) \cap \R$, for all $\theta \in
\Gamma-\R_+$. Moreover, given a non-threshold eigenvalue
$\lambda_0 \in \sigma(\Delta)$, the eigenspace
$E_{\lambda_0}(\Delta)$, associated to $\Delta$ and $\lambda_0$,
has finite dimension bounded by the degree of the pole $\lambda_0$
of the map $\lambda \mapsto R(\lambda,\theta)$. This algebraic
multiplicity is independent of $\theta \in \Gamma-\R_+$.
\\
\\
b) Fix $\theta \in \Gamma$. For $f,g \in \mathscr{V}$ the function
$$\lambda \mapsto \langle R(\lambda)f,g \rangle_{L^2(Z,E)}$$
has a meromorphic continuation from $\Lambda$ to $\C-\big(
\sigma_{ess}(\Delta_\theta)\cup \sigma_{d}(\Delta_\theta)\bigl)$,
where $\sigma_{ess}(\Delta_\theta)$  is the set computed in
section \ref{section:essential spectrum cyl end}.
\\
\\
c)  $\Delta$ has no singular spectrum.
\\
\\
d) Let $\theta_1, \theta_2 \in \Gamma$ be such that $arg(\theta'_1) \geq arg(\theta'_0)$ for $0<arg(\theta'_i)< \pi/2 $, we have:
\begin{equation}
\sigma_{d}(\Delta_{\theta_0})=\sigma_{d}(\Delta_{\theta_1}) \cap \sigma_{d}(\Delta_{\theta_0}).
\end{equation}
e) Non-thresholds eigenvalues of $\Delta$ are isolated (with
respect to the eigenvalues of $\Delta$) and may only accumulate on
the set of thresholds\footnote{We prove in corollary
\ref{cor:infty unique accum point cyl end} that, in fact, $\infty$
is the unique possible accumulation point.} or at $\infty$.
\\
\\
f) If the lowest eigenvalue, $\mu_0$, of $\Delta_Y$ is larger than $0$ then $\sigma_d(\Delta)$ is a discrete subset of $[0,\mu_0)$. The unique possible accumulation point of $\sigma_d(\Delta)$ is $\gamma_0$. If $\mu_0=0$, then $\sigma_d(\Delta)=\emptyset$; in other words, all eigenvalues are embedded in the continuous spectrum.
\end{thm}
The next proposition  follows from the definition of essential
spectrum and (\ref{eq:essspec1}).
\begin{prop}
i) If $\lambda \in \sigma_{pp}(\Delta_\theta)$ and $\lambda \notin
\sigma_{ess}(\Delta_\theta)$, then $\lambda $ is an isolated
eigenvalue of finite multiplicity.
\\
\\
ii) For $Re(\theta)>0$, $\sigma_{pp}(\Delta_\theta)$ accumulates in $\sigma_{ess}(\Delta_\theta)$. In particular, the real part of the pure point spectrum of $\Delta_\theta$ accumulates only in $\sigma(\Delta_Y)$.
\end{prop}
Next we show that the unique possible accumulation point of
$\sigma_{pp}(\Delta)$ is $\infty$. For that we use the following
theorem (see \cite{DONELLY2},pag. 352).
\begin{thm}\cite{DONELLY2} \label{thm: Donelly Weyl Law cyl. end} If $N(\lambda)$ denotes the number of eigenvalues of $\Delta$ which are less than $\lambda$, then one has
\begin{equation}
N(\lambda)\leq C \lambda^{m-1/2}.
\end{equation}
\end{thm}
In \cite{DONELLY2} the previous theorem is proved only for the
Laplacian $\Delta$ acting on functions, but  it generalizes easily
to our context. As we have previously said, theorem  \ref{thm:
Donelly Weyl Law cyl. end} implies the following corollary.
\begin{cor}\label{cor:infty unique accum point cyl end}
The unique possible accumulation point of $\sigma_{pp}(\Delta)$ is $\infty$.
\end{cor}
We define the set of resonances of $\Delta$ in the following way:
\begin{equation}\label{eq: def resonances1}
\mathscr{R}_\theta(\Delta)=\{ \lambda \in
\sigma_{d}(\Delta_\theta): \lambda \notin \sigma_{pp}(\Delta)\}.
\end{equation}
The parameter $\theta$ is simply uncovering new pure point
spectrum in the sense of the following proposition, that is a
consequence of the uniqueness of the meromorphic extension of
$\lambda \mapsto \left< R(\lambda)f,g \right>_{L^2(Z,E)}$, for
$f,g \in \mathscr{V}$.
\begin{prop}
Suppose $\theta_1, \theta_0 \in \Gamma$ and
$0<arg(\theta_0)<arg(\theta_1)<\pi/2$. Then
\begin{equation}
\mathscr{R}_{\theta_0}(\Delta)\subset
\mathscr{R}_{\theta_1}(\Delta).
\end{equation}
\end{prop}
There is an analogue version of the previous proposition for
$\theta_1, \theta_0 \in \Gamma$ and
$0>arg(\theta_0)>arg(\theta_1)>-\pi/2$.
\\
\\
We recall some facts about the analytic  extension of the
resolvent $R(\lambda)$ of $\Delta$. First, we introduce some
notation. Let $\Sigma$ be the Riemann surface on which the
functions $\sqrt{z-\mu_i}$ are defined. Observe that $\Sigma$ is a
$\omega$-covering of $\C$, with ramification points $\{\mu_i:i \in
\N\}$. Denote:
\begin{equation}
\begin{split}
L^2&_\delta(Z,E):=\\
&\{\varphi:Z \to E: \text{ measurable section s.t. } \int_0^\infty \int_Y \normv{\varphi}^2e^{2 \delta u}dvol(y)du < \infty\}.
\end{split}
\end{equation}
In \cite{GUILLOPE} (see also \cite{HUS}) the resolvent is extended
as a function of $\lambda \in \Sigma$ taking values in the bounded
operators from $L^2_{-\delta}(Z,E)$ to $L^2_{\delta}(Z,E)$. The
next theorem provides more information about the  resonances of
$\Delta$:
\begin{thm}(\cite{HUS}, theorem 3.26)
Suppose that  $\lambda \in \mathscr{R}_\theta(\Delta)$. Then:
\begin{itemize}
\item[1)] Suppose $Im(\lambda) \neq 0$. Then, if $\lambda \notin
\bigcup_{i \in \N}\mu_i+\theta'[0,\infty)$ then $\lambda \in
\sigma_d(\Delta_\theta)$. Under these conditions, if
$0<arg(\theta')<\pi/2$, then $Im(\lambda)>0$; if
$0>arg(\theta')>-\pi/2$, then $Im(\lambda)<0$. \item[2)] There are
not real resonances different  than the set of thresholds
$\tau(\Delta_\theta)=\sigma(\Delta_Y)$. If $\lambda=\mu_i$ for
some $i \in \N$, then the resolvent, as a function from
$L_\delta^2(Z,E)$ to $L_{-\delta} ^2(Z,E)$, has a pole of at most
second order. In fact, $\mu_i$ is a pole of second order always
that it is a $L^2$-eigenvalue of $\Delta$; in this case, the
leading part of the Laurent expansion of $R(\lambda)$ at $\mu_i$
is the orthogonal projection in the $L^2$-eigenspace space
$E_{\mu_i}$.
\end{itemize}
\end{thm}
As a consequence of the previous theorem, we have the following
corollary that completes the description of the spectrum of
$\Delta$ of theorem  \ref{thm: descriptio delta theta cyl}.
\begin{cor}\label{cor: descriptio delta theta 2}
The real resonances of $\Delta$ are contained in $\sigma(\Delta_Y)$.
\end{cor}

\subsection{$\Delta_\theta$ are $m$-sectorial}
\label{sec: Deltatheta is sectorial}The  theory of $m$-sectorial
operators and forms that we use in this section is described in
appendix \ref{chap:Ichinose}. Our goal  is to show that the
operators $\Delta_\theta$, for $\theta \in \Gamma$ (see (\ref{def:
Gamma})), are $m$-sectorial (see definition \ref{def: strictly
sectorial}). This result will be important when we calculate
$\sigma_{ess}(H_\theta)$.
\\
\\
Let $\eta_0 \in C^\infty(\R_+)$ be a positive real function such
that $\eta_0(u)=1$ for $u<1$, $\eta_0(u)=0$ for $u \in
[K-1,\infty)$, and $\eta'_0(u)\leq 0$ for $u\in[1,K-1]$, where we
are considering $K>2$. Let $\eta_1:=1-\eta_0$. Both $\eta_0$ and
$\eta_1$ induce functions on $Y \times \R_+$, defining
$\eta_k(u,y):=\eta_k(u)$ for $k=0,1$; making $\eta_k$ equal to $0$
where it is not defined, we can extend it to all of $Z$. In this
way we think $\eta_0$ and $\eta_1$ as functions in $C^\infty(Z)$.
\begin{prop} \label{prop:pruebaichinoselema2}
For  $Re(\theta)>0 $ there exist a $\gamma \geq 0$  such that
\begin{equation}
\begin{split}
 Re \left( \langle a_0( \theta, u) f\right. & \left.,\eta_1 f \rangle_{L^2(Z,E)}+\gamma \langle f,f\rangle _{L^2(Z,E)}\right)\geq\\
& \vert Im(\langle a_0( \theta, u) f,\eta_1 f\rangle_{L^2(Z,E)})
\vert
\end{split}
\end{equation}
for all $f \in L^2(Z,E)$.
\end{prop}
{\bf Proof:}
\\
We observe that, by definition, $\eta_0+\eta_1=1$. Let $\gamma
\geq 0$, then, for all $z \in Z$:
\begin{equation}
\begin{split}
&Re\left(\langle a_0( \theta, u) f,\eta_1 f\rangle(z)+\gamma \langle f,f\rangle (z)\right)-\vert Im(\langle a_0( \theta, u) f,\eta_1 f \rangle (z)) \vert  \\
&\hspace{1cm}\geq \eta_1(z) \langle f,f \rangle (z) \{ Re(a_0(\theta,u))-\vert Im(a_0( \theta, u))\vert+\gamma\}+\gamma \eta_0 \langle f,f\rangle (z).
\end{split}
\end{equation}
Notice  that in the previous calculations the inner product denote
the Hermitian product in the fiber $E_z$. Since, for all $z \in Z$
$\eta_1(z) \langle f,f\rangle (z)$ and $\gamma \eta_0(z)\langle
f,f \rangle (z)$ are both equal or larger than $0$, then it is
enough to prove that there exist $\gamma>0$ such that
\begin{equation}
\begin{split}
Re(a_0(\theta,u)-\normv{Im(a_0( \theta, u))}+\gamma \geq 0
\end{split}
\end{equation}
for all $u \in \R_+$. This is true because $\{a(\theta,u): u \in \R_+\}$ is a compact subset of $\R^2$ (for $\theta$ fixed), any compact subset of $\R^2$ is inside a cube $[-n,n]^2$, and we can always find a $N$ such that $[N-n,n+N]^2$ is inside a cone, with slope $1$, included in a right-half-plane. $\square$
\\
\\
Let us now turn back to the set $\Gamma$  defined in (\ref{def:
Gamma}). The next theorem shows that the operators $\Delta_\theta$
are $m$-sectorial for $\theta \in \Gamma$. We will use this fact
when we compute the essential spectrum of $H_\theta$,  in section
\ref{section:calculating essent. spec}.
\begin{thm} \label{thm:sectorial} For $\theta \in \Gamma$
there exists a $\gamma(\theta) \in \R_+$ such that the form with domain  $\mathscr{W}_1(Z,E)$ defined by  $f \mapsto \langle \Delta_\theta f,f \rangle _{L^2(Z,E)}+\gamma(\theta) \langle f,f\rangle _{L^2(Z,E)} $ is  $m$-sectorial.
\end{thm}
{\bf Proof:}
\\
Let us prove that there exist $k>0$ and $\gamma \in \R$ such that
for all $f \in \mathscr{W}_2(Z,E)$:
\begin{equation} \label{ineq:sectorial}
Re\left(\langle \Delta_\theta f,f\rangle_{L^2(Z,E)}+\gamma \langle
f,f \rangle _{L^2(Z,E)} \right)\geq  k\vert Im\left(\langle
\Delta_\theta f,f\rangle_{L^2(Z,E)}\right) \vert.
\end{equation}
We observe that having the previous inequality,  the theorems
\ref{thm:sectop1} and \ref{thm:sectop2} imply that the bilinear
form $f \mapsto \langle \Delta_\theta f,f\rangle_{L^2(Z,E)}+\gamma
\langle f,f\rangle_{L^2(Z,E)} $ is strictly $m$-sectorial, and
hence the form defined by $\Delta_\theta$ is $m$-sectorial.
\\
\\
Next we prove the inequality (\ref{ineq:sectorial}).  We use
proposition \ref{prop: formula htheta cylend} to see that
\begin{equation}\label{eq:msectorial1}
\begin{split}
\langle  \eta_0& \Delta_\theta (f),f\rangle _{L^2(Z,E)}=\langle
\parcial{}{u}(f), \parcial{}{u}(\eta_0)f\rangle
_{L^2(Z,E)}+\langle \nabla (f),\eta_0 \nabla (f)\rangle _{L^2(Z,E
\otimes T^ *Z)},
\end{split}
\end{equation}
for $f \in \mathscr{W}_1(Z,E)$. Now let $a_2(\theta,
u)\parcial{^2}{u^2}+a_1(\theta,u)\parcial{}{u}+a_0(\theta,u)+\Delta_Y$
be the local expression of the operator $\Delta_\theta$ (see
equation  (\ref{eq:a1 a2 a3 theta})). We have:
\begin{equation}\label{eq:msectorial2}
\begin{split}
&\langle \eta_1 \Delta_\theta (f),f \rangle _{L^2(Z,E)}=-\langle \parcial{}{u}(f),\parcial{}{u}\left(\overline{a}_2
\eta_1\right) f\rangle _{L^2(Z,E)}-\langle \parcial{}{u}(f),\left( \overline{a}_2\eta_1 \right)\parcial{}{u}(f)\rangle _{L^2(Z,E)}\\
&+\langle \parcial{}{u}(f),\overline{a}_1\eta_1 f\rangle_{L^2(Z,E)}+\langle f,\overline{a}_0(u)\eta_1 f\rangle_{L^2(Z,E)}+\langle \eta_1 \Delta_Y (f),f\rangle_{L^2(Z,E)}.
\end{split}
\end{equation}
Using  (\ref{eq:msectorial1}) and (\ref{eq:msectorial2}) we find
\begin{equation}
\begin{split}
&\langle \eta_0 \Delta_\theta (f),f\rangle _{L^2(Z,E)}+\langle\eta_1 \Delta_\theta (f),f\rangle_{L^2(Z,E)}\\
&=\langle\parcial{}{u}(f), \left(-\parcial{}{u}(\eta_0)-\parcial{}{u}\left(\overline{a}_2 \eta_1\right)+\overline{a}_1 \eta_1 \right)f\rangle_{L^2(Z,E)}\\
&\hspace{0.8cm}-\langle \parcial{}{u}(f),\left(\overline{a}_2 \eta_1 \right)\parcial{}{u}(f)\rangle_{L^2(Z,E)}+\langle \nabla (f),\eta_0 \nabla (f)\rangle_{L^2(Z,E)}\\
&\hspace{0.8cm}+\langle f,\overline{a}_0 \eta_1 f\rangle+\langle
\eta_1\Delta_Y (f),f\rangle_{L^2(Z,E)}.
\end{split}
\end{equation}
Since
\begin{equation}
\begin{split}
\langle \nabla (f), \eta_0 \nabla (f)\rangle_{L^2(Z,E)}=&\int_{Z_0}\langle \nabla (f), \eta_0 \nabla (f)\rangle(z)dz\\
&+\int_{Y\times [0,\infty)}\langle \nabla_Y (f), \eta_0 \nabla_Y (f)\rangle(z)dz\\
&-\int_{Y\times [0,\infty)}\langle \parcial{}{u} (f),  \eta_0
\parcial{}{u}  (f)\rangle dz,
\end{split}
\end{equation}
we have that the term
\begin{equation}
\begin{split}
s(f):=&\int_{Z_0}\langle \nabla (f), \eta_0 \nabla (f)\rangle (z)dz+\int_{Y\times [0,\infty)}\langle \nabla_Y (f), \eta_0 \nabla_Y (f)\rangle (z)dz\\
&+\langle \eta_1\Delta_Y (f),(f)\rangle _{L^2(Z,E)}
\end{split}
\end{equation}
is greater or equal than $0$.
\\
\\
We define the bilinear form $h(\theta)$ by
\begin{equation}
\begin{split}
h(\theta)(f):&=\langle \parcial{}{u}(f), \left(-\parcial{}{u}(\eta_0)-\parcial{}{u}\left(\overline{a}_2(\theta) \eta_1\right)+\overline{a}_1(\theta) \eta_1\right)f\rangle_{L^2(Z,E)}\\
&\hspace{0.5cm}-\langle \parcial{}{u}(f),\left( \overline{a}_2(\theta) \eta_1 \right)\parcial{}{u}(f)\rangle_{L^2(Z,E)}\\
&\hspace{0.5cm}+\int_{Y\times [0,\infty)}\langle \parcial{}{u} (f),  \eta_0 \parcial{}{u}  (f)\rangle(z)dz.
\end{split}
\end{equation}
Since $\langle \Delta_\theta
(f),f\rangle_{L^2(Z,E)}=h(\theta)(f)+s(f)+\langle f,\overline{a}_0
\eta_1 f\rangle_{L^2(Z,E)}, $  $s(f) \geq 0$, theorem \ref{thm:sum
of sectorials}, proposition \ref{prop:pruebaichinoselema2} and the
definition of $h(\theta)$, to finish the proof of
(\ref{ineq:sectorial}), it only remains to prove that there exist
$\gamma>0$ and $k>0$ that satisfy
\begin{equation}\label{ineq: htheta is sectorial}
\begin{split}
Re(h(\theta)f)+\gamma \langle f,f\rangle _{L^2(Z,E)} \geq k\normv{Im (h(\theta)f)}.
\end{split}
\end{equation}
We observe that
$-\parcial{}{u}(\eta_0)-\parcial{}{u}\left(\overline{a}_2(\theta)\eta_1\right)+\overline{a}_1(\theta)\eta_1$
has support on $[0,\infty)$ (as a function of $u$, $\theta$ fixed)
and it is bounded there by a constant $C$ (see proposition
\ref{prop:aibounded}). Then:
\begin{equation*}
\begin{split}
Re(h(\theta)(f))&-k \vert Im(h(\theta)(f)) \vert  \geq \\
&\int_{Y\times [0,\infty)}\{(-Re(a_2(\theta))+k \vert Im(a_2(\theta)) \vert) \eta_1 +\eta_0\}\langle \parcial{}{u}f,\parcial{}{u}f)\rangle(z)dz\\
&-C\int_{Y \times [0,\infty]} \langle \vert \parcial{}{u} f \vert,\vert f \vert \rangle (z)dz.
\end{split}
\end{equation*}
Notice that $Re(a_2)\eta_1=1/2Re(a_2)\eta_1
+1/2Re(a_2)-1/2Re(a_2)\eta_0$. Using the fact that, for all
$\epsilon \in \R$,
$$
\epsilon^2 \vert u \vert^2+1/4\epsilon^2 \vert v\vert^2\geq \vert u \vert  \vert v \vert,
$$
we have for all $k \in \R_+$:
\begin{equation} \label{eq: mother msectorial}
\begin{split}
Re(h(\theta)(f))&-k \vert Im(h(\theta)(f)) \vert  \geq \\
&\int_{Y\times [0,\infty)}(-1/2Re(a_2(\theta))+k \vert Im(a_2(\theta)) \vert) \eta_1 \langle \parcial{}{u}(f),
\parcial{}{u}(f)\rangle (z)dz\\
&+ \int_{Y \times [0,\infty]}\bigl(-1/2Re(a_2(u,\theta))-C \epsilon^2\bigr) \langle \parcial{}{u}(f),
\parcial{}{u} (f) \rangle (z)dz\\
&+ \int_{Y \times [0,\infty]} \bigl(1+1/2Re(a_2(u,\theta))\bigr)\eta_0\langle \parcial{}{u}(f),\parcial{}{u} (f) \rangle (z)dz\\
&-C/(4\epsilon^2) \int_{Y \times [0,\infty]} \langle f , f \rangle (z)dz.
\end{split}
\end{equation}
Recall that $a_2(\theta,u)=\frac{-1}{\psi'_\theta(u)^2}$. We
observe that $\psi'_\theta$ is bounded and
$\psi'_\theta(u)=B(u)\theta+1$ where
$B(u):=\varphi'(u)u+\varphi(u) \geq 0$ for all $u \in [0,\infty)$.
For $\theta \in \Gamma$, $Im(\theta)^2<1/2$, and
$(B(u)Re(\theta)+1)^2 \geq 1$, then:
\begin{equation}
(B(u)Re(\theta)+1)^2-Im(\theta)^2 \geq 1-Im(\theta)^2>1/2.
\end{equation}
Thus, there exists a constant $C_0>0$ such that
\begin{equation}
\begin{split}
-Re(a_2(u,\theta))&=\frac{(B(u)Re(\theta)+1)^2-Im(\theta)^2}{\vert
\psi'_\theta(u)\vert^4} \\
&\geq \frac{1}{2\max \{u\in
[0,\infty):\vert \psi'_\theta(u) \vert\}^4}>C_0>0
\end{split}
\end{equation}
for all $u \in \R_+$. Hence,  we can find $\epsilon $ such that
\begin{equation}\label{eq: calculations for mother msectorial1}
(-1/2Re(a_2(\theta,u))-C \epsilon^2 )>0.
\end{equation}
Now we show that, for all $\theta \in \Gamma$, there exists  $k$
such that $Re(\psi_\theta'(u)^2)-k \normv{Im(\psi_\theta'(u)^2} \geq 0$, for all $u \in \R_+$.
We observe that $\psi_\theta'(u)^2=B(u)^2\theta^2+2B(u)\theta+1$.
Suppose that $\theta:=\theta_0+i\theta_1$, for $\theta_0$ and $\theta_1$ real numbers. We denote by
\begin{equation}
M:=\max\{B(u):u\in \R_+\}=\max\{\varphi(u)u+\varphi'(u):u\in \R_+\}.
\end{equation}
From the definition of $\Gamma$ in (\ref{def: Gamma}) it follows that, for $\theta \in \Gamma$, $\theta_0^2- \normv{\theta_1}^2\geq 0$, then:
\begin{equation}
\begin{split}
Re(&\psi_\theta'(u)^2)-k \normv{Im(\psi_\theta'(u)^2)}\geq 1-k\{2M^2\normv{\theta_0\theta_1}+2M\normv{\theta_1}\}.
\end{split}
\end{equation}
The previous calculations show that, for fixed $\theta \in \Gamma$, we can always find a $k$ such that
\begin{equation}\label{eq: calculations for mother msectorial2}
Re(\psi_\theta'(u)^2)-k \normv{Im(\psi_\theta'(u)^2} \geq 0,
\end{equation}
for all $u \in \R_+$. Finally, we observe that:
\begin{equation}\label{eq: calculations for mother msectorial3}
1+1/2Re(a_2(\theta,u))=1-1/2\frac{Re(\overline{\psi'_\theta}^2(u))}{\vert \psi'_\theta(u)\vert^4} \geq 0,
\end{equation}
because $2\vert \psi'_\theta(u)\vert^4-Re(\overline{\psi'_\theta}^2(u)) \geq 0$.
This last
inequality is true because, for all $a \in \C$, $2\normv{a}^2 \geq Re( \overline{a})$.
\\
\\
Using  (\ref{eq: mother msectorial}), (\ref{eq:
calculations for mother msectorial1}), (\ref{eq: calculations for
mother msectorial2}) and (\ref{eq: calculations for mother
msectorial3}), we can conclude that there exists $K_0 \geq 0$ such that
\begin{equation}
\begin{split}
Re(h(\theta)(f))&-k \vert Im(h(\theta)(f)) \vert  \geq K_0-C/(4\epsilon^2) \int_{Y \times [0,\infty]} \langle f , f \rangle(z)dz.
\end{split}
\end{equation}
Finally, we can take $\gamma$ large enough to have:
\begin{equation}
\begin{split}
Re(h(\theta)(f))&-k \vert Im(h(\theta)(f)) \vert+\gamma\langle f,f\rangle_{L^2(Z,E)}  \geq 0,
\end{split}
\end{equation}
what finishes the proof of (\ref{ineq: htheta is sectorial}), and
with it the proof of the theorem.$\square$

\section{Analytic dilation on complete manifolds with corners of codimension 2}
\label{subsec:Manifolds with corners}Let $X$ be a complete
manifold with corners of codimension 2 and let $E$ be an Hermitian
vector bundle over $X$.  Let $\Delta$ be a  generalized Laplacian
acting on $C^\infty(X,E)$. We say that $\Delta$ is a {\bf compatible
Laplacian} over $X$ if it satisfies the following properties:
\begin{itemize}
\item[i)] There exist Hermitian vector bundle $E_i$ over $Z_i$
such that $E \vert_{[0,\infty)\times Z_i}$ is the pull-back of
$E_i$ ($i=1,2$). We suppose also that the Hermitian metric of $E$
is the pullback of the Hermitian metric of $E_i$.  In addition, on $[0,\infty)
\times Z_i$ we have
\begin{equation}
\Delta=-\parcial{^2}{u_k^2}+\Delta_{Z_i},
\end{equation}
where $\Delta_{Z_i}$ is a compatible  Laplacian acting on
$C^\infty(Z_i,E_i)$. \item[ii)] There exists  Hermitian vector
bundle $S$ over $Y$ such that $E \vert_{[0,\infty)^2 \times Y}$ is
the pull-back of $S$, and on $[0,\infty)^2 \times Y$ we have,
\begin{equation}
\Delta=-\parcial{^2}{u_1^2}-\parcial{^2}{u_2^2}+\Delta_Y
\end{equation}
where $\Delta_Y$ is a generalized Laplacian acting on $C^\infty(Y,S)$.
\end{itemize}
Since $X$ is a  manifold with bounded geometry and the vector
bundle $E$ has bounded Hermitian metric and bounded connection,
$\Delta:C^\infty_c(X,E) \to L^2(X,E)$ is essentially self-adjoint
(see \cite{Shubin}).  We denote by $H$ its self-adjoint extension.
For $i=1,2$, $\Delta_{Z_i}:C^\infty_c(Z_i,E_i) \to L^2(Z_i,E_i)$
is also essentially self-adjoint and we denote its self adjoint
extension by   $H^{(i)}$. Let $b_i$ be the self-adjoint extension
of $-\parcial{^2}{u_i^2}:C^\infty_c([0,\infty)) \to
L^2([0,\infty))$ obtained with Dirichlet boundary conditions. We
denote $H_i$ the self-adjoint operator $-b_i\otimes 1+1 \otimes
H^{(i)}$ acting on $L^2(\R)\otimes L^2(Z_i,E_i)$. Similarly
$H^{(3)}$ denotes the self-adjoint operator associated to the
essentially self-adjoint operator $\Delta_Y:C^\infty_c(Y,S) \to
L^2(Y,S)$, and we denote by $H_3$, the self-adjoint operator
$H_3:=-b_1\otimes 1\otimes 1-1\otimes b_2\otimes 1+ 1\otimes
1\otimes H^{(3)}$ acting on $L^2([0,\infty)) \otimes L^2([0,\infty))
\otimes L^2(Y)$.
\\
\\
The operators $H_i$ are called  {\bf channel operators} for
$i=1,2$ and $3$. The operators $H_1$ and $H_2$ have a free
channel of dimension 1 (associated to $b_1$ and $b_2$
respectively), $H_3$ is channel operator with a free channel of
dimension 2 (associated to $-b_1\otimes 1\otimes 1-1\otimes
b_2\otimes 1$).  In some parts of this text we abuse of the
notation and denote by $H$ the  Laplacian acting on distributions.
\subsection{The definition of $U_\theta$ for $\theta \in [0,\infty)$}\label{def:Utheta corner}
For $i=1,2$ and $\theta \in [0,\infty)$,
$U_{i,\theta}:L^2(Z_i,E_i) \to L^2(Z_i,E_i)$ denotes an analytic
dilation operator associated to the Laplacian $H^{(i)}$. The unitary operator
$U_{i,\theta}$ was  described in section \ref{analytic dilation
cyl end}. In this section we denote $H^{(i),\theta}:=U_{i,\theta}
H^{(i)} U_{i,\theta}^{-1}$.
\\
\\
In the next definition we use the exhaustion defined in
(\ref{eq:def exhaust}). Let $\theta \in [0,\infty)$:
\begin{equation}
\begin{split}
U_\theta f(x):=\begin{cases}f(x_0)\hspace{1.5cm} \text{ for } x=x_0 \in X_0,\\
\\
f(m_i, \psi_\theta(u_i))\psi_\theta'^{1/2}(u_i)\\
\hspace{2.6cm} \text{ for } x=(m_i,u_i) \in M_i \times [0,\infty), i=1,2,\\
\\
\\
f(y,\psi_\theta(u_1),\psi_\theta(u_2))\psi_\theta'^{1/2}(u_1)\psi_\theta'^{1/2}(u_2)\\
\hspace{2.6cm}\text{for  } x=(y,u_1,u_2) \in Y\times [0,\infty)
\times [0,\infty).
\end{cases}
\end{split}
\end{equation}
The following proposition follows from the definition of
$U_\theta$.
\begin{prop}Let $\theta \in [0,\infty)$ and $f \in C^\infty_c(X,E)$,
\begin{itemize}
\item[i)]  $U_\theta f \in C^\infty_c(X,E)$.
\item[ii)] $U_\theta$ extends to a unitary operator in $L^2(X,E)$.
\end{itemize}
\end{prop}
The inverse of $U_\theta$ is given by:
\begin{equation}
\begin{split}
U_\theta^{-1}f(x):=\begin{cases} f(x_0) \text{ for } x=x_0 \in X_0,\\
\\
f(m_i, \alpha_\theta(u_i))\psi_\theta'^{-1/2}(\alpha_\theta(u_i))\\
 \hspace{2.6cm} \text{ for }x=(m_i,u_i) \in M_i \times [0,\infty), i=1,2.\\
\\
f(y,\alpha_\theta(u_1),\alpha_\theta(u_2))\psi_\theta'^{-1/2}(\alpha_\theta(u_1))\psi_\theta'^{-1/2}(\alpha_\theta(u_2))\\
\hspace{2.6cm} \text{  for } x=(y,u_1,u_2) \in Y\times
[0,\infty)^2.
\end{cases}
\end{split}
\end{equation}
One can check the following proposition.
\begin{prop}\label{prop:rotinzi}For $\theta \in [0,\infty)$ and $i=1,2$, if $f \in C^\infty_c(X,E)$ and $(z_i,u_i) \in  [0,\infty) \times Z_i$, then:
\\
i) $U_\theta f(u_i,z_i)=(U_{i,\theta}f)(\psi_\theta(u_i),z_i)\psi_\theta'^{1/2}(u_i)$.
\\
\\
ii) $U_\theta^{-1} f(u_i,z_i)=U_{i,\theta}^{-1}f(\alpha_\theta(u_i),z_i)\psi_\theta'^{-1/2}(\alpha_\theta(u_i)).$
\end{prop}
\subsection{The family $H_\theta$ for $\theta \in \Gamma$}\label{sec:famila Htheta}
For $\theta \in [0,\infty)$, we define the operator
$H_\theta:=U_\theta H U_\theta^{-1}$. By direct calculation, one
can prove:
\begin{prop}\label{prop: deltatheta1}
For $\theta \in [0,\infty)$ and $i=1,2$, if $f \in
C^\infty_c(X,E)$, then for all $(z_i,u)\in Z_i \times [0,\infty)$
\begin{equation}\label{eq:defdeltatheta}
\begin{split}
H_\theta f(u,z_i)=&H^{(i),\theta} f(u,z_i)-(\parcial{^2}{u^2}f))(u,z_i)\alpha_\theta'(\psi_\theta(u))^2\\
&-(\parcial{}{u}f)(u,z_i)\alpha_\theta''(\psi_\theta(u))\\
&+(\parcial{}{u}f)(u,z_i)\psi_\theta'(u)^{-1}\psi_\theta''(u)\alpha_\theta'(\psi_\theta(u))^2\\
&-3/4f(u,z_i)\psi_\theta'(u)^{-3/2}\psi_\theta''(u)^2\alpha_\theta'(\psi_\theta(u))^2\\
&+1/2f(u,z_i)\psi_\theta'(u)^{-1}\psi_\theta'''(u)\alpha_\theta'(\psi_\theta(u))^2\\
&+1/2f(u,z_i)\psi_\theta'(\alpha_\theta(u))^{-1}\psi_\theta''(u)\alpha_\theta''(\psi_\theta(u)).
\end{split}
\end{equation}
In particular, if $u_i>R$, we have:
\begin{equation}
H_\theta
f(u_i,z_i)=-\theta'\parcial{^2}{u_i^2}f(u_i,z_i)+H^{(i),\theta}f(u_i,z_i),
\end{equation}
and, if $u_i<K$,
\begin{equation}
H_\theta f(u_i,z_i)=-\parcial{^2}{u_i^2}f(u_i,z_i)+H^{(i)}f(u_i,z_i).
\end{equation}
\end{prop}
Similarly, the next proposition describe the operator $H_\theta$
on $Y \times [0,\infty)^2$:
\begin{prop}\label{propextdeltatheta2}
For $\theta \in [0,\infty)$, for all $f \in C^\infty_c(X,E)$ and
for all $(y,u_1,u_2)\in Y \times [0,\infty)^2$,
\begin{equation}
\begin{split}
H_\theta f(y,u_1,u_2)=&H^{(3)} f(y,u_1,u_2))-\sum_{i=1}^2\{(\parcial{^2}{u_i^2}f))(y,u_1,u_2)\alpha_\theta'(\psi_\theta(u_i))^2\\
&-(\parcial{}{u}f)(y,u_1,u_2)\alpha_\theta''(\psi_\theta(u_i))\\
&+(\parcial{}{u}f)(y,u_1,u_2)\psi_\theta'(u_i)^{-1}\psi_\theta''(u_i)\alpha_\theta'(\psi_\theta(u_i))^2\\
&-3/4f(y,u_1,u_2)\psi_\theta'(u_i)^{-3/2}\psi_\theta''(u_i)^2\alpha_\theta'(\psi_\theta(u_i))^2\\
&+1/2f(y,u_1,u_2)\psi_\theta'(u_i)^{-1}\psi_\theta'''(u_i)\alpha_\theta'(\psi_\theta(u_i))^2\\
&+1/2f(y,u_1,u_2)\psi_\theta'(\alpha_\theta(u_i))^{-1}\psi_\theta''(u_i)\alpha_\theta''(\psi_\theta(u_i))\}.
\end{split}
\end{equation}
In particular, if $u_1, u_2>R$, we have:
\begin{equation}
H_\theta
f(y,u_1,u_2)=-\theta'\parcial{^2}{u_1^2}f(y,u_1,u_2)-\theta'\parcial{^2}{u_2^2}f(y,u_1,u_2)+H^{(3)}f(y,u_1,u_2).
\end{equation}
an if $u_1, u_2<K$ we have:
\begin{equation}
H_\theta
f(y,u_1,u_2)=-\parcial{^2}{u_1^2}f(y,u_1,u_2)-\parcial{^2}{u_2^2}f(y,u_1,u_2)+H^{(3)}f(y,u_1,u_2).
\end{equation}
\end{prop}
We observe that, for all $f \in C^\infty_c(X,E)$ and $(y,u_1,u_2)
\in Y \times [0,\infty)^2$, we can write:
\begin{equation}
\begin{split}
H_\theta (f)(y,u_1,u_2)=&\sum_{i=1}^2 \{a_2(\theta,u_i)\parcial{^2}{u_i^2}f(y,u_1,u_2)\\
&+a_1(\theta,u_i)\parcial{}{u_i}f(y,u_1,u_2)+a_0(\theta,u_i)f(y,u_1,u_2)\}
\end{split}
\end{equation}
where the functions $(\theta,u) \mapsto a_i(\theta,u)$ were defined
in (\ref{eq:a1 a2 a3 theta}).
\\
\\
Our next goal is to prove that $(H_\theta)_{\theta \in \Gamma}$ is
an holomorphic family of type A. We use the general theory of
uniformly elliptic operators on manifolds with bounded geometry, as
explained for example in  \cite{Shubin}.
\begin{prop}\label{prop:Htheta elliptic}
For all $\theta \in \Gamma$, the operator $H_\theta$ is an uniformly elliptic operator, that is,
\begin{equation}
\normv{\sigma_2(H_\theta)^{-1}(x,\xi)}^{-1} \leq \normv{\xi}^{-2},
\end{equation}
for all $x \in X$ and $\xi \in T_xX$.
\end{prop}
As a consequence of  proposition \ref{prop:Htheta elliptic} we have the next corollary.
\begin{cor}\label{cor:equivnormX}
The operators $H_\theta$ are closed operators in the domain $\mathscr{W}_2(X,E)$.
\end{cor}
Given $f \in Dom(H)$ and $g \in L^2(X,E)$, we will
prove that the function $\theta \mapsto \langle H_\theta f,g
\rangle_{L^2(X,E)}$ is holomorphic for $\theta \in \Gamma$, accordingly we consider  the following partition of unity of $X$. Let
$\eta \in C^\infty_c(\R)$ be such that $\eta(u)=1$ for $u \leq
K-2$ and $\eta(u)=0$ for $u \geq K-1$. Let $\kappa:=1-\eta$, we
define the following functions with their natural extensions to
the whole $X$. Let $(z_i,u_i) \in Z_i \times [0,\infty)$, then:
\begin{equation*}
\begin{split}
\kappa_i(z_i,u_i):=\kappa(u_i), \hspace{0.5cm} \eta_i(z_i,u_i):=1-\kappa_i.
\end{split}
\end{equation*}
We observe that $\eta_i+\kappa_i=1$. In particular, we have that
$(\eta_1+\kappa_1)(\eta_2+\kappa_2)=1$. We study the functions
$\theta \mapsto \kappa_1 \kappa_2 H_\theta$ and $\theta \mapsto
\eta_i \kappa_j H_\theta$. The next proposition is a technical
tool.
\begin{prop}
1) Given $f \in Dom(H)$ and $g \in L^2(X,E)$, there exists $h \in L^1(X)$ such that:
\begin{equation}
\normv{< \kappa_1 \kappa_2 H_\theta f,g>(x) }\leq h(x),
\end{equation}
and
\begin{equation} \label{eq:deltathetalebesguethm}
\normv{\parcial{}{\theta}(< \kappa_1 \kappa_2 H_\theta f,g>(x))} \leq h(x),
\end{equation}
for $\theta$ in any compact subset of $\Gamma$.
\\
\\
2) For $i,j \in \{1,2\}$, $i \neq j$, and $\theta$ in a compact
subset of $\Gamma$ given $f \in Dom(H)$ and $g \in L^2(X,E)$,
there exists $h \in L^1(X)$ such that:
\begin{equation}
\normv{< \kappa_i \eta_j H_\theta f,g> (x)}\leq h(x),
\end{equation}
and
\begin{equation} \label{eq:deltathetalebesguethm2}
\normv{\frac{d}{d\theta}(<\kappa_i \eta_j H_\theta f,g>(x))} \leq h(x).
\end{equation}
\end{prop}
As a consequence of the previous proposition we can apply well known
theorems (see \cite{BAUER}, page 89) to put partial
derivatives inside of the respective integrals. This and the
Cauchy-Riemann equations imply the following theorem.
\begin{thm}\label{thm: H holom}
Given $f \in Dom(H)$ and $g \in L^2(X,E)$, the function $\theta \mapsto \langle H_\theta f ,g \rangle_{L^2(X,E)}$ is holomorphic.
\end{thm}
Observe that theorem \ref{thm: H holom} and corollary \ref{cor:equivnormX} prove theorem \ref{thm:Afamily} in the introduction.
\subsection{The essential spectrum of $H_\theta$}
\label{section:calculating essent. spec}
The goal of this section is to prove theorem \ref{thm: ess Htheta} in the introduction. Let us consider the  set:
\begin{equation}
\begin{split}
\mathscr{F}_\theta:=& \left(\bigcup_{i=1}^2 \bigcup_{\lambda \in \sigma_{pp}(H^{(i)})} \lambda +\theta'[0,\infty) \right)\\
&\cup \left( \bigcup_{\mu \in \sigma(H^{(3)})} \mu +\theta'[0,\infty) \right)\\
& \cup \left(\bigcup_{i=1}^2 \bigcup_{\gamma \in
\mathscr{R}(H^{(i),\theta})} \gamma +\theta'[0,\infty) \right),
\end{split}
\end{equation}
where $\mathscr{R}(H^{(i),\theta})$ is defined by:
\begin{equation}\label{eq: def resonances corner}
\mathscr{R}(H^{(i),\theta}):=\{\lambda \in
\sigma_{pp}(H^{(i),\theta}): \lambda \notin \sigma_{pp}(H^{(i)})
\}.
\end{equation}
The elements of $\mathscr{R}(H^{(i),\theta})$ shall be called {\bf
resonances of $H^{(i),\theta}$}. We observe that the set $\mathscr{R}(H^{(i),\theta})$
is independent of $\theta$ in the sense that if $arg(\theta_1')
\geq arg(\theta_2')$ then $\sigma_{pp}(H^{(i),\theta'_2})\subset
\sigma_{pp}(H^{(i),\theta'_1})$ (see  item d), theorem  \ref{thm: descriptio
delta theta cyl}).  The next proposition follows easily
from the definition of $\mathscr{F}_\theta$.
\begin{prop}
For $\theta \in \Gamma$, the following equation holds:
\begin{equation}\label{eq:prop Ftheta=Gtheta}
\begin{split}
\mathscr{F}_\theta=& \left(\bigcup_{i=1}^2 \bigcup_{\lambda \in \sigma_{pp}(H^{(i),\theta})} \lambda +\theta'[0,\infty)\right)\\
&\cup \left( \bigcup_{\mu \in \sigma(H^{(3)})} \mu
+\theta'[0,\infty) \right).
\end{split}
\end{equation}
\end{prop}
By the previous proposition theorem \ref{thm: ess Htheta} is equivalent to $\sigma_{ess}(H_\theta)=\mathscr{F}_\theta$ for $\theta \in \Gamma$. We use
the results of appendix \ref{chap:DHSV}, about the sets
$N_\infty(A)$ (see definition  \ref{defin: Ninfty}), $N_{ess}(A)$
(see definition \ref{def: Ness}) and $\sigma_{ess}(A)$ to show that $\sigma_{ess}(H_\theta)=\mathscr{F}_\theta$. The proof has
two basic steps:
\begin{itemize}
\item[i)] To prove that $N_{\infty}(H_\theta)=\mathscr{F}_\theta$.
\item[ii)] To prove that there exists $\eta_0^d \in C^\infty_c(X)$ such that $supp (\eta_0^d) \subset X_d$, and for all $f \in C^\infty_c(X,E)$,
\[
\norm{[H_\theta,\eta_0^d]f}_{L^2(X,E)}\leq \epsilon(d)(\norm{H_\theta f}_{L^2(X,E)}+\norm{f}_{L^2(X,E)})
\]
with $\epsilon(d)\to 0$ as $d \to \infty$.
\end{itemize}
Using ii) we see that $H_\theta$ satisfies the conditions for
applying theorem \ref{thm1:DHSV},  and we get
$N_{ess}(H_\theta)=N_\infty(H_\theta)$. Using i) and part iii) of
theorem  \ref{thm:DHSVNess} in appendix \ref{chap:DHSV}, we prove
$\sigma_{ess}(H_\theta)=\mathscr{F}_\theta$.
\subsubsection{The equality $N_{\infty}(H_\theta)=\mathscr{F}_\theta$}
The boundary Weyl sequences (abbreviately bWs), defined in definition
\ref{defin: Ninfty}, will play a very important role in this
section . Let $\mu \in \sigma(H^{(3)})$ and $\lambda \in
\mu+\theta'[0,\infty)$. We observe that we can apply theorems
\ref{thm:DHSVNess} and \ref{thm1:DHSV} to the operators
$\theta'\frac{d^2}{du^2}$ and $\theta'\frac{d^2}{du^2}+\mu$. Then
there exist  a bWs, $(p_n)$,  associated to $0$ and the operator
$\theta'\parcial{^2}{u_1^2}$,  and a bWs, $(q_n)$,  associated to
$\lambda$ and the operator $\theta' \parcial{^2}{u_2^2}+\mu$. Let
$\varphi \in C^\infty(Y,S)$ be a normal eigenfunction of $H^{(3)}$
associated to the eigenvalue $\mu$.
\begin{prop}\label{prop:glueingbWs1}
Let $g_n\in C^\infty_c( Y\times [0,\infty)^2,E)$  be  defined by
$g_n(u_1,u_2,y)=p_n(u_1)q_n(u_2)\varphi(y)$. Then $(g_n)$ induces
a boundary Weyl sequence for $\lambda$ and the operator
$H_\theta$.
\end{prop}
{\bf Proof:}
\\
It is easy to check $g_n \in C_c^\infty(X,E)$, $\norm{g_n}=1$ and
that, for all $K>0$, there exists a $N$ such that, for all $n>N$,
$supp (g_n) \cap X_K =\emptyset$, where $X_K$ is defined in
(\ref{eq:def exhaust}). We observe that:
\begin{equation}\label{eq: pn qn varphi is bws}
\begin{split}
&\norm{\left(\theta'\parcial{^2}{u_1^2}+\theta'\parcial{^2}{u_2^2}+H^{(3)}-\lambda \right)g_n} \leq \\
&C(\norm{\theta'\parcial{^2}{u_1^2}p_n}+\norm{(\theta'\parcial{^2}{u_2^2}+\mu-\lambda) q_n}).
\end{split}
\end{equation}
Since $(p_n)$ is a bWs of $\theta'\parcial{^2}{u_1^2}$ and the value
$0$, and $(q_n)$ is a bWs of $-\theta'\parcial{^2}{u_2^2}+\mu$ and
the value $\lambda$, the last two terms of  (\ref{eq: pn qn varphi
is bws}) tend to $0$. We have proved that $\lim_{n \to
\infty}\norm{\left(H_\theta-\lambda \right)g_n}=0$. $\square$
\\
\\
Now let $\gamma \in \sigma_{pp}(H^{(i),\theta})$ and $\lambda \in
\gamma + \theta'[0,\infty)$.  Let $\varphi \in C^\infty(Z_i,E_i)$
be a normal $L^2$-eigenfunction of $H^{(i),\theta}$ with
eigenvalue $\gamma$. Let $\eta \in C^\infty(\R)$ such that
$\eta(u)=1$, for $u\leq 1$; $\eta(u)=0$, for $u>2$; and
$\eta'(u)\leq 0$. Define $\eta_n(u):=\eta(\frac{u}{n})$. Let $f_n$
be a bWs associated to the operator $-\theta'\frac{d^2}{du^2}$ and
$0$.
\begin{prop} \label{prop:glueingbWs2}We denote\footnote{ With this notation $u_j$ is the
real variable in the cylinder of $Z_i$.} $i,j \in\{1,2\}$ such
that $i \neq j$. Let $g_n \in C^\infty(Z_i \times [0,\infty),E)$
be defined by
$g_n(u_1,u_2,z_i):=\frac{1}{\norm{\eta_n(u_j)f_n(u_i)\varphi(z_1)}}\eta_n(u_j)f_n(u_i)\varphi(z_i)$.
Then, $(g_n)$ induces a boundary Weyl sequence associated to
$H_\theta$ and the value $\lambda$.
\end{prop}
{\bf Proof:}
\\
It is easy to check $g_n \in C^\infty_c(X,E)$, $\norm{g_n}=1$, and that for all $K \in \N$ there exists an $N$ such that for all $n \geq N$ $supp g_n \cap X_K=\emptyset$. Since $\lim_{n\to \infty}\frac{1}{\norm{\eta_n(u_2)f_n(u_1)\varphi(z_1)}_{L^2(X,E)}}=1$, it is enough to prove:
$$\lim_{n \to \infty}\norm{(\theta'\parcial{^2}{u_i^2}+H^{(i),\theta})g_n}=0.$$
We observe that:
\begin{equation}\label{eq:gluing from R+XZi to X}
\begin{split}
(H^{(i),\theta}-\gamma)(\eta_n(u_j)\varphi(z_i))=&-\theta' \parcial{^2}{u_j^2}(\eta_n)\varphi-2\theta'\parcial{}{u_j}(\eta_n)\parcial{}{u_j}(\varphi).
\end{split}
\end{equation}
Since $\eta_n(u_j)=\eta(\frac{u_j}{n})$, then
$\parcial{}{u_j}(\eta_n)(u_j)=\frac{1}{n}\parcial{}{u_j}(\eta)(\frac{u_j}{n})$
and
$\parcial{^2}{u_j^2}(\eta_n)(u_j)=\frac{1}{n^2}\parcial{^2}{u_j^2}(\eta)(\frac{u_j}{n})$.
Hence, using equation  (\ref{eq:gluing from R+XZi to X})
$$\norm{(\theta'\parcial{^2}{u_i^2}+H^{(i),\theta}-\gamma)g_n}\leq C(A_n+B_n+C_n),$$
where
\begin{equation}
\begin{split}
&A_n:=\norm{\eta_n(u_j)\varphi(z_i)(\theta'\parcial{^2}{u_i^2})f_n(u_i)};\\
&B_n:=\norm{\left(\theta'\parcial{^2}{u_j^2}(\eta_n)(u_j)\right)\varphi(z_i)f_n(u_i)}=\norm{\left(\theta'\parcial{^2}{u_j^2}(\eta_n)(u_j)\right) \varphi(z_i)};\\
&C_n:=\norm{\left(\theta'\parcial{}{u_j}(\eta_n)(u_j)\right)\varphi(z_i)f_n(u_i)} =\norm{\left(\theta'\parcial{}{u_j}(\eta_n)(u_j)\right) \varphi(z_i)}.
\end{split}
\end{equation}
For $A_n$, we have:
\begin{equation}
\begin{split}
A_n \leq  \norm{(\theta'\parcial{^2}{u_i^2})f_n(u_i)}\to 0.
\end{split}
\end{equation}
For $B_n$, we estimate:
\begin{equation}
\begin{split}
B_n&\leq C(\theta) \frac{1}{n^2} \left(\int_{Z_i}\vert (\parcial{^2}{u_j^2}\eta)(\frac{u_j}{n})\varphi(z_i)\vert^2 dz_i\right)^{1/2}\leq \frac{1}{n^2}\norm{\varphi} \to 0;
\end{split}
\end{equation}
and, finally, for $C_n$:
\begin{equation}
\begin{split}
C_n\leq C(\theta)
 \frac{1}{n} \left(\int_{Z_i}\vert (\parcial{}{u_j}\eta)
 (\frac{u_j}{n})\varphi(z_i)\vert^2 dz_i\right)^{1/2} \leq \frac{1}{n} \norm{\varphi} \to 0. \text{  } \square
\end{split}
\end{equation}
Propositions \ref{prop:glueingbWs1} and \ref{prop:glueingbWs2}
prove that $\mathscr{F}_\theta \subset N_\infty(H_\theta)$. Now we
are going to prove the other inclusion. We recall that we  denote by
$b_i$ the self-adjoint operator in $L^2([0,\infty))$ obtained from
$-\parcial{^2}{u_i^2}$ with Dirichlet boundary conditions.  We
denote by $H_{i,\theta}$  the closed operator $-\theta' 1\otimes
b_i+H^{(i),\theta}\otimes 1$ acting on $L^2(Z_i \times
[0,\infty),E)=L^2(Z_i,E_i) \otimes L^2([0,\infty))$.
\begin{prop}
For $i=1,2$, and $\theta \in \Gamma$:
\\
i) $\sigma_{ess}(H_{i,\theta})=N_\infty(H_{i,\theta})$.
\\
ii) $N_\infty(H_{i,\theta})$ is given by:
\begin{equation}\label{eq:Ninfty Hi theta}
\begin{split}
N_\infty(H_{i,\theta})=\left(\bigcup_{\lambda \in
\sigma_{pp}(H^{(i),\theta})} \lambda +\theta'[0,\infty)
\right)\cup \left( \bigcup_{\mu \in \sigma (H^{(3)})}  \mu
+\theta'[0,\infty)\right).
\end{split}
\end{equation}
\end{prop}
{\bf Proof:}
\\
In the same way that we proved $\mathscr{F}_\theta \subset
N_\infty(H_\theta)$,  using propositions \ref{prop:glueingbWs1}
and \ref{prop:glueingbWs2},  we can prove
\begin{equation}\label{eq:Ftheta subset Ninfty Hi theta}
\left(\bigcup_{\lambda \in \sigma_{pp}(H^{(i),\theta})}  \lambda
+\theta'[0,\infty)  \right) \cup \left( \bigcup_{\mu \in \sigma
(H^3)}  \mu +\theta'[0,\infty) \right) \subset
N_\infty(H_{i,\theta}).
\end{equation}
We denote by $\mathscr{F}_{i,\theta}$  the right-side of
(\ref{eq:Ninfty Hi theta}). Then, (\ref{eq:Ftheta subset Ninfty Hi
theta}) implies $\mathscr{F}_{i,\theta} \subset
N_\infty(H_{i,\theta})$, and by proposition \ref{prop:N subset
Ninfty}, $\mathscr{F}_{i,\theta} \subset
N_\infty(H_{i,\theta})\subset \sigma_{ess}(H_{i,\theta})$. We know
that the operator $H^{(i),\theta}$ is $m$-sectorial (theorem
\ref{thm:sectorial}), so we can apply Ichinose lemma (see theorem
\ref{thm:Ichinoselemma}) in the next computations:
\begin{equation}
\begin{split}
\sigma(H_{i,\theta})&=\sigma(H^{(i),\theta})+\sigma(-\theta' \parcial{^2}{u_i^2})\\
&=\left(\theta' [0,\infty) +
\sigma_{pp}(H^{(i),\theta})\right)\cup \left( \sigma(H^{(3)})
+\theta' [0,\infty)\right).
\end{split}
\end{equation}
The above equation implies $N_\infty(H_{i,\theta})\subset
\mathscr{F}_{i,\theta}$ and $\sigma_{ess}(H_{i,\theta})\subset
\mathscr{F}_{i,\theta}$, that together with  (\ref{eq:Ftheta
subset Ninfty Hi theta}) implies
$N_\infty(H_{i,\theta})=\mathscr{F}_{i,\theta}=\sigma_{ess}(H_{i,\theta})$.
We have proven the proposition. $\square$
\\
\\
Next  we prove $N_{\infty}(H_\theta) \subset \mathscr{F}_\theta$.
Let $\lambda \in N_{\infty}(H_\theta)$ and let $f_n \in
C^\infty_c(X,E)$ be a bWs associated to the operator $H_\theta$
and $\lambda$. In order to prove $N_\infty(H_\theta) \subset
\mathscr{F}_\theta $ we will construct, using $f_n$, a bWs
associated to $\lambda$ and one of the operators $H_{i,\theta}$.
We define $\kappa_n:=1-\eta_n$.
\begin{prop}
There exists $c>0$ and a subsequence $s$ of $\N$ such that
\begin{equation}
\norm{\kappa_n(u_1)f_{s(n)}}_{L^2(X,E)}\geq c>0 \hspace{0.25cm}\text{ or } \hspace{0.25cm} \norm{\kappa_n(u_2)f_{s(n)}}_{L^2(X,E)}\geq c>0.
\end{equation}
\end{prop}
{\bf Proof:}
\\
Suppose $\norm{\kappa_n(u_1)f_n}_{L^2(X,E)}\to 0$ and $\norm{\kappa_n(u_2)f_n}_{L^2(X,E)}\to 0$. We can choose $s$ such that $\chi_{n}f_{s(n)}=f_{s(n)}$ where $\chi_n$ denotes the characteristic function of $X-X_{n^2+1}$. Since $\chi_n^2\leq \kappa_n(u_1)^2+\kappa_n(u_1)^2$, then
\begin{equation}
\begin{split}
&1=\norm{\chi_{n}f_{s(n)}}_{L^2(X,E)}^2\leq  \norm{\kappa_n(u_1)f_{s(n)}}_{L^2(X,E)}^2+\norm{\kappa_n(u_2)f_{s(n)}}_{L^2(X,E)}^2,
\end{split}
\end{equation}
which is a contradiction. $\square$
\\
\\
The previous proposition allows us to suppose that $0<c<\norm{\kappa_n(u_1)f_{s(n)}}_{L^2(X,E)}$.
\begin{prop}Denote by $g_n$ the function in $C^\infty(X,E)$, defined by
$g_n:=\frac{1}{\norm{\kappa_n(u_1)f_{s(n)}}_{L^2(X,E)}}\kappa_n(u_1)f_{s(n)}$.
Then, $g_n$ induces  a boundary Weyl sequence associated to
$H_{1,\theta}$ and $\lambda$.
\end{prop}
{\bf Proof:}
\\
It is easy to check that $\norm{g_n}_{L^2(X,E)}=1$  and, for all $T>0$,
there exists $N \in \R$ such that $\forall n \geq N$, $supp g_n
\cap X_T=\emptyset$. Denoting $\kappa(u):=1-\eta(u)$, we define
$\kappa_n(u_1):=\kappa(\frac{u_1}{n})$, then
$\parcial{}{u_1}(\kappa_n)(u_1)=\frac{1}{n}\parcial{}{u_1}(\kappa)(\frac{u_1}{n})$
and
$\parcial{^2}{u_1^2}(\kappa_n)(u_1)=\frac{1}{n^2}\parcial{^2}{u_1^2}(\kappa)(\frac{u_1}{n})$.
Hence:
\begin{equation}
\begin{split}
\norm{\left(H_{1,\theta}-\lambda\right)(\kappa_n(u_1)f_{s(n)})}_{L^2(X,E)}\leq A_n +B_n+C_n \to 0,
\end{split}
\end{equation}
where
\begin{equation*}
\begin{split}
A_n^2&:=4\norm{\parcial{}{u_1}(\kappa_n)(u_1)\parcial{}{u_1}(f_{s(n)})}_{L^2(X,E)}^2\\
&\hspace{0.4cm} \leq \frac{C}{n^2} \int_{Z_1 \times [0,\infty)}\normv{\parcial{}{u_1}(\kappa)(\frac{u_1}{n})\parcial{}{u_1}(f_{s(n)})}^2 dvol(x)\\
&\hspace{0.4cm} \leq \frac{C}{n^2} \norm{ \parcial{}{u_1}
(f_{s(n)}) }^2 \leq C \norm{f_{s(n)}}_{2}^2 \cdot \frac{1}{n^2}
\leq C' \frac{1}{n^2} \to 0.
\end{split}
\end{equation*}
In the last inequalities we use that $\norm{ \parcial{}{u_1}
(f_{s(n)}) } \leq  C \norm{f_{s(n)}}_{2}$, which follows from the
theory of bounded differential operators on manifolds with bounded
geometry (see  \cite{Shubin}, \cite{EICHHORN}). We use also
$\norm{f_n}_2 \leq C$, which follows from $\norm{f_n}=1$, $\lim_{n
\to \infty}\norm{(H_\theta -\lambda)f_n} = 0$ and since $H_\theta$
is uniformly elliptic. For $B_n$, we have:
\begin{equation*}
\begin{split}
B_n^2:=\norm{\parcial{^2}{u_1^2}(\kappa_n)(u_1)f_{s(n)}}^2 \leq & C \int_{Z_1 \times [0,\infty)}\normv{\parcial{^2}{u_1^2}(\kappa_n)(u_1)f_{s(n)}}^2 dvol(x)\\
& \leq C \frac{1}{n^4} \norm{f_{s(n)}}=C \frac{1}{n^4} \to 0.
\end{split}
\end{equation*}
Finally for $C_n$, we have
\begin{equation*}
\begin{split}
C_n:=\norm{\kappa_n(u_1)\left(H_{1,\theta}-\lambda\right)f_{s(n)}}
\leq C \norm{\left(H_{\theta}-\lambda \right)f_{s(n)}} \to
0.\text{  }\square
\end{split}
\end{equation*}
Next we prove step ii) of our proof of $\sigma_{ess}(H_\theta)=\mathscr{F}_\theta$.  Let $\eta \in C^\infty(\R)$ such that
$\eta(u)=1$ for  $u\leq 1$, $\eta(u)=0$ for $u>2$ and
$\eta'(u)\leq 0$. Denote $\eta_n(u):=\eta(\frac{1}{n})$,  and
define
\begin{equation}
\eta_0^{(d)}(u_1,u_2,y):=\eta_d(u_1)\eta_d(u_2).
\end{equation}
\begin{prop}
For all $f \in C^\infty_c(X,E)$,
\begin{equation}
\norm{[H_\theta,\eta_0^{(d)}]f}_{L^2(X,E)}\leq
\epsilon(d)(\norm{H_\theta f}_{L^2(X,E)}+\norm{f}_{L^2(X,E)}),
\end{equation}
with $\lim_{d\to \infty} \epsilon(d)=0$.
\end{prop}
{\bf Proof:}
\\
Let $i,j \in \{1,2\}$ and $i \neq j$. We observe that
\begin{equation}
\begin{split}
H_\theta (\eta_0^{(d)} f)=&\theta'\sum_{i=1}^1 \{ 2\eta_d(u_j)\parcial{}{u_i}(\eta_d)(u_i)\parcial{}{u_i}(f)\\
&+\eta_d(u_j)\parcial{^2}{u_i^2}(\eta_d)(u_i)f\}\\
&+\eta_d(u_j)\eta_d(u_i)H_\theta(f).
\end{split}
\end{equation}
Hence,
\begin{equation}
\begin{split}
\norm{[H_\theta,\eta_0^{(d)}]f}_{L^2(X,E)} \leq
 &\normv{\theta'}\cdot\sum_{i=1}^1 \{ 2\norm{\eta_d(u_j)\parcial{}{u_i}(\eta_d)(u_i)\parcial{}{u_i}(f)}_{L^2(X,E)}\\
&+\norm{\eta_d(u_j)\parcial{^2}{u_i^2}(\eta_d)(u_i)f}_{L^2(X,E)}.
\end{split}
\end{equation}
By definition of $\eta_0^{(d)}$,
$\parcial{}{u_1}(\eta_0^{(d)})(u_1,u_2,y)=\frac{1}{d}\parcial{}{u_1}(\eta)(u_1/d)\eta_d(u_2)$.
Thus,
\begin{equation}
\begin{split}
&\norm{\eta_d(u_j)\parcial{}{u_i}(\eta_d(u_i))\parcial{}{u_i}(f)}_{L^2(X,E)}^2 \leq  \int_X \vert \eta_d(u_j)\parcial{}{u_i}(\eta_d)(u_i)\parcial{}{u_i}(f) \vert^2 dx\\
&\hspace{2cm}\leq 1/d^2 \int_{Y \times {[0,\infty)^2}} \vert
\parcial{}{u_1}(\eta)(u_1/d)\eta_d(u_2) \parcial{}{u_i}(f)\vert^2
dx\\
&\hspace{2cm} \leq
\frac{C}{d^2}(\norm{f}_{L^2(X,E)}+\norm{H_\theta f}_{L^2(X,E)})^2
\to 0.
\end{split}
\end{equation}
In the last inequality we use that
$\parcial{}{u_1}(\eta)(u_1/d)\eta_d(u_2) \parcial{}{u_i}(f)$ is a
bounded differential operator of degree 1 (hence a continuous
operator from $\mathscr{W}_2(X,E)$ to $L^2(X,E)$), and the fact
that the norm $f \mapsto \norm{f}_{L^2(X,E)}+\norm{H_\theta
f}_{L^2(X,E)}$ is equivalent to $\norm{\cdot}_2$ from the theory
of uniformly elliptic operators acting on sections of vector
bundles on manifolds with bounded geometry (see \cite{Shubin}). We
finish the proof of the proposition with the following
calculation:
\begin{equation}
\begin{split}
\norm{\eta_d(u_j)\parcial {^2}{u_i^2}(\eta_d)(u_i)f}_{L^2(X,E)}^2
\leq 1/d^4 \int_{Y \times {[0,\infty)^2}} \vert f \vert^2 dx \to
0.\text{ }\square
\end{split}
\end{equation}
\subsection{Analytic vectors}\label{sect:analytic vect corners}
In this section we construct a  subset $\mathscr{V}$ of $L^2(X,E)$
such that
\begin{itemize}
\item[i)] $\mathscr{V}$ is a dense subset of $L^2(X,E)$.
\item[ii)] For $ f \in \mathscr{V}$ the function $\theta \mapsto U_\theta f \in L^2(X,E)$ is well defined for $\theta \in \C$, $Re(\theta)>0$.
\item[iii)] $U_\theta \mathscr{V}$ is dense in $L^2(X,E)$ for $Re(\theta)>0$.
\end{itemize}
Let $\mathscr{V}_i$ be the analytic vectors associated to
$U_{i,\theta}$ (see equation (\ref{eq:analytvectors1})). Let $\eta
\in C^\infty([0,\infty))$ be such that $\eta'\geq 0$ and
$$
\eta(u)=\begin{cases}
1&K<u<\infty.\\
0&0<u\leq K-1.
\end{cases}
$$
For $i=1,2,$ define $\eta_i(z_i,u_i):=\eta(u_i)$ and extend them to  $X$. Denote $\kappa:=1-\eta_1-\eta_2$. Define
\begin{equation}
\mathscr{V}:=\{(\kappa g+\sum_{i=1}^2\eta_i \frac{1}{u_i^2}p_i(\frac{1}{u_i})f_i(z_i):g\in L^2(X,E),p_i(x)\in \C[x] \text{ and } f_i \in \mathscr{V}_i\}.
\end{equation}
Then i), ii) and iii) can be deduced from the properties of  $\mathscr{V}_i$ given in section \ref{analytic vectors cyl end}.
\subsection{Consequences of Aguilar-Balslev-Combes theory}
In section \ref{section:calculating essent. spec}, we calculated
the essential spectrum of $H_\theta$. The  following theorem is a
consequence of theorem  \ref{thm: essential spectrum of Htheta}
and the existence of the analytic vectors satisfying i), ii) and
iii) of section \ref{sect:analytic vect corners}. It can be
proved using the general ideas of  Aguilar-Balslev-Combes theory
as explained in \cite{HISLOPSIGAL}.
\begin{thm}\label{thm: descriptio H theta cor}We have:
\\
\\
a) The  set  of non-threshold\footnote{The set of thresholds of
$H$, $\tau(H)$, is equal to $\sigma(H^3) \cup \bigcup_{i=1}^2
\sigma_{pp}(H^{(i),\theta})$.} eigenvalues of $H$ is equal to
$\sigma_d(H_\theta) \cap \R$, for all $\theta \in
\Gamma-[0,\infty)$. Moreover, given a non-threshold
eigenvalue $\lambda_0$, the eigenspace $E_{\lambda_0}(H)$, associated to $H$
and $\lambda_0$, has finite dimension bounded by the degree of the
pole $\lambda_0$ of the map $\lambda \mapsto R(\lambda,\theta)$.
This algebraic multiplicity is independent of $\theta \in
\Gamma-[0,\infty)$.
\\
\\
b) Fix $\theta \in \Gamma$. For $f,g \in \mathscr{V}$ the function
$$\lambda \mapsto \langle R(\lambda)f,g\rangle_{L^2(X,E)}$$
has a meromorphic continuation from $\Lambda$ to
$\C-\left(\sigma_{ess}(H_\theta)\cup\sigma_{pp}(H_\theta)\right)$,
where $\sigma_{ess}(H_\theta)$ was calculated in theorem \ref{thm:
essential spectrum of Htheta}.
\\
\\
c)  $H$ has no singular spectrum.
\\
\\
d) Let $\theta_1, \theta_2 \in \Gamma$ be such that $arg(\theta'_1) \geq arg(\theta'_0)$ for $0<arg(\theta'_i)< \pi/2 $, we have:
\begin{equation}
\sigma_{d}(H_{\theta_0})=\sigma_{d}(H_{\theta_1}) \cap \sigma_{d}(H_{\theta_0}).
\end{equation}
e) Non-thresholds eigenvalues of $H$ are isolated (with respect to the eigenvalues of $H$) and may only accumulate on the set of thresholds or at $\infty$.
\\
\\
f) If the lowest threshold, $\gamma_0$, is larger than $0$, then $\sigma_d(H)$ is a discrete subset of $[0,\gamma_0)$. In this case, the unique possible accumulation point of $\sigma_d(H)$ is $\gamma_0$. If $\gamma_0=0$, then $\sigma_d(H)= \emptyset$, in other words all eigenvalues are embedded in the continuous spectrum.
\end{thm}
At the moment we do not know if there is a compatible  Laplacian
having embedded eigenvalues. It seems that the natural
conjecture is that generically there are no embedded eigenvalues.
Similarly, we do not know if it is possible to find a generalized
Laplacian that has embedded eigenvalues accumulating at one of the
thresholds. We believe that it is possible to prove that
$\sigma_{pp}(H_\theta)$ can accumulate only from below  in
$\tau(H)$, and we will address this question in a forthcoming paper.  In
particular, this would imply that $0$ is not an accumulation point
of $\sigma_{pp}(H_\theta)$.
\\
\\
The next proposition  is easy to prove from the definition of
essential spectrum and from theorem \ref{thm: essential spectrum
of Htheta}.
\begin{prop}(\cite{CANOTHESIS})
i) If $\lambda \in \sigma_{pp}(H_\theta)$ and $\lambda \notin \sigma_{ess}(H_\theta)$, then $\lambda $ an isolated eigenvalue of finite multiplicity.
\\
\\
ii) For $Re(\theta)>0$, the accumulation points of $\sigma_{pp}(H_\theta)$ are contained in  $\sigma_{ess}(H_\theta)$. In particular, the real part of the pure point spectrum of $H_\theta$ can accumulate only in $\tau(H)$.
\end{prop}
In \cite{CANOGENEIGEN}, using the results of this section, we
define generalized eigenfunctions associated to
$L^2$-eigenfunctions of $H^{(1)}$, $H^{(2)}$ and $H^{(3)}$. The
wave packets associated to these generalized eigenfunctions
describe completely the absolutely continuous spectrum,
$L^2_{ac}(X,E)$, associated to $H$. In fact, in
\cite{CANOASYMPCOMPLETENESS} is proven the asymptotic completeness
of waves operators associated to $H_1$, $H_2$ and $H_3$; the
generalized eigenfunctions of \cite{CANOGENEIGEN} describe such
wave operators.

\appendix
\section{Geometric spectral analysis of $\sigma_{ess}$} \label{chap:DHSV}
In this appendix we remark that the results of section 3 of the
paper \cite{DHSV} are also valid in the context of generalized
Laplacians on complete manifolds with corners of codimension 2.
In fact  the proofs of  the theorems in this appendix are
essentially the same that those given in  \cite{DHSV}, and because
of that we refer readers to that paper or to \cite{CANOTHESIS}.
\\
\\
We begin by recalling the definition of singular sequences associated to a closed operator $A$:
\begin{defin}\label{def:singseq}
A sequence $(f_n)_{n \in \N}\subset Dom(A)$ is {\bf a singular sequence for $A$ associated to the value $\lambda \in \C$} if and only if
\\
\\
i) $\norm{f_n}=1$ and $(f_n)_{n \in \N}$ has no $L^2$-convergent subsequence.
\\
\\
ii) $\lim_{n \to \infty}\norm{(A-\lambda)f_n}=0$.
\end{defin}
In this section we distinguish between different types of singular
sequences of a geometric operator and we describe some relations
between them. They define different subsets of the essential
spectrum defined in  (\ref{eq: def ess d pp}).
\\
\\
Let $A:Dom(A) \subset L^2(X)\to L^2(X)$ be a closed operator.
\begin{defin}\label{def: Ness}
Define the set {\bf $N_{ess}(A)$} of $\lambda \in \C$ such that there exists a sequence $(u_n)_{n \in \N} \subset Dom(A)$ such that $\norm{u_n}=1$, $u_n \to 0$ (weakly) and $\norm{(\lambda-A)u_n}\to 0$.
\end{defin}
Observe that if $\lambda \in N_{ess}(A)$, then the sequence
$(u_n)_{n \in \N}$ associated to $\lambda$ is a singular sequence
in the sense of definition  (\ref{def:singseq}).
\\
\\
Now we define another important class of singular sequences $N_\infty(A)$.
\begin{defin}(\cite{DHSV}, page 10)\label{defin: Ninfty}
Let {\bf $N_\infty(A)$} be the set of $\lambda \in \C$ such that there exists a sequence $(u_n)_{n \in \N} \subset C^\infty_c(X)$ with
\begin{itemize}
\item[i)] $\norm{u_n}=1$,
\item[ii)] $\norm{(A -\lambda)u_n}\to 0$,
\item[iii)] for every  compact subset  $K \subset M$ there exists  $N \in \N$ such that for $n >N$, $supp (u_n) \cap K = \emptyset$.
\end{itemize}
We will call the sequence $u_n$ a {\bf boundary Weyl sequence} (abbr. bWs).
\end{defin}
Observe that if  $\lambda \in N_\infty(A)$ and  $(u_n)_{n \in \N}$
is a sequence as in definition \ref{defin: Ninfty}, then $(u_n)_{n
\in \N}$ is a singular sequence associated to $A$ and the value
$\lambda$.
\begin{prop}(\cite{DHSV}, page 9)\label{prop:N subset Ninfty}
i) $N_{ess}(A) \subset \sigma_{ess}(A)$.
\\
\\
ii) $N_{ess}(A)$ is closed.
\end{prop}
The following theorem gives conditions for the equality of $\sigma_{ess}(A)$ and $N_{ess}(A)$.
\begin{thm}\label{thm: Weyl criterion bws}(\cite{DHSV}, theorem 3.1)(Weyl's criterion for $\sigma_{ess}(A)$)
\label{thm:DHSVNess} Let $A$ be a closed operator on a Hilbert space $\mathscr{H}$ with non-empty resolvent set. Then:
\begin{itemize}
\item[i)]$ N_{ess}(A) \subset\sigma_{ess}(A)$.
\item[ii)] The boundary of $\sigma_{ess}(A)$ is contained in $N_{ess}(A)$.
\item[iii)] $N_{ess}(A)=\sigma_{ess}(A)$ if and only if each connected component of the complement of $N_{ess}(A)$ contains a point of $\rho(A)$.
\end{itemize}
\end{thm}
The next theorem gives conditions for the equality of
$N_{\infty}(A)$ and $N_{ess}(A)$. We will use the notation $X_0$
and $X_d$ for the manifolds defined in  (\ref{eq:def exhaust}) for
$T=0$ and $T=d$ respectively.
\begin{thm}\label{thm1:DHSV}(\cite{DHSV}, theorem 3.2)
Let $A$ be a closed operator on $L^2(X)$ with non-empty resolvent set, having $C^\infty_c(X)$ as a core. Let $\eta_0 \in C^\infty_c(X)$ such that $\eta_0(x)=1$ for $x \in X_0$. Let $\eta_0^d \in C^\infty_c(X)$ such that $\eta^d_0(x)=1$ for $x \in X_d$, and $0 \leq \eta^d_0(x) \leq 1$ for $x \in X$. Suppose $\forall d$, $\eta_0^d(z-A)^{-1}$ is compact for some $z \in \rho(A)$, and that for all $u \in C^\infty_c(X)$,
\begin{equation}\label{eq: hyp estim commutator}
\norm{[A,\eta_0^d]u}\leq \epsilon(d)(\norm{Au}+\norm{u}),
\end{equation}
with $\epsilon(d)\to 0$ as $d\to \infty$. Then $N_\infty(A)=N_{ess}(A).$
\end{thm}
\section{Ichinose lemma}\label{chap:Ichinose}
In this appendix we recall some definitions and we formulate the Ichinose lemma (theorem \ref{thm:Ichinoselemma}).
The following definitions follow \cite{RS1} and \cite{KATO}, we refer there for a deeper study of the topic. Let $q$ be a bilinear form on a Hilbert space $\mathscr{H}$ with domain $Q(q)$.
\begin{defin}\label{def:closed form}
We say that $q$ is {\bf closed} if and only if always that a sequence $\varphi_n \in Q(q)$ converges   $\varphi$ in the norm topology, and
$$lim_{n,m\to \infty}q(\varphi_n-\varphi_m,\varphi_n-\varphi_m)=0,$$
then, we have $\varphi \in Q(q)$ and $q(\varphi_n-\varphi,\varphi_n-\varphi)\to 0$.
\end{defin}
\begin{defin}\label{defin:sect}
 A quadratic form $q$ is {\bf sectorial} if there exists  $\theta$, $0<\theta<\pi/2$ with $\vert arg(q(\varphi,\varphi))\vert \leq \theta$ for all $\varphi \in Q(q)$.
\end{defin}
\begin{defin}
A quadratic form $q$ is called {\bf strictly $m$-accretive} if it is both closed and sectorial.
\end{defin}
\begin{defin}\label{def: strictly sectorial}
A form $q$ is called  {\bf strictly $m$-sectorial} if there are complex numbers $z$ and $e^{i\alpha}$, with $\alpha$ real, so that $e^{i\alpha}q+z$ is strictly $m$-accretive. The operator $T$ associated to $q$ is also called {\bf strictly $m$-sectorial}.
\end{defin}
Observe that in order to prove that $q$ is strictly $m$-sectorial it is enough to show that there exists $\gamma \in \R$ and $k \in \R_+$ such that for all $f \in Q(q)$
\begin{equation}
k Re(q(f))- \vert Im(q(f)) \vert \geq \gamma (f,f).
\end{equation}
Every closed operator $T$ defines a dense form $q(T)$ by
\begin{equation}\label{eq: form from operator}
q(t)(\varphi,\psi):=(\varphi,T\psi),
\end{equation}
for $\varphi, \psi \in D(T)$.
\begin{defin}
An operator $T$ is sectorial if there is a $\theta$, $0<\theta<\pi/2$ such that its numerical range, $\Theta(T)$, is a subset of a sector   $\{z \in \C: \vert arg(z)\vert \leq \theta\}$.
\end{defin}
The following theorems are important in section \ref{sec: Deltatheta is sectorial}.
\begin{thm}(\cite{KATO}, page 318)\label{thm:sectop1}
A sectorial operator $T$ is form closable, that is, the form
$q(T)$ defined by (\ref{eq: form from operator}) has an extension
that is closed in the sense of definition \ref{def:closed form}.
\end{thm}
\begin{thm}(\cite{KATO},page 316)\label{thm:sectop2}
Let $\tilde{q}$ be the closure of a densely defined form $q$. The numerical range $\Theta(q)$ of $q$ is a dense subset of the numerical range $\Theta(\tilde{q})$ of $\tilde{q}$.
\end{thm}
The next theorem is also used in section \ref{sec: Deltatheta is sectorial}.
\begin{thm}(\cite{KATO}, page 319)\label{thm:sum of sectorials}
Let $q_1,...q_s$ be sectorial forms in $\mathscr{H}$ and let $q:=q_1+...+q_s$ [with $D(q):=D(q_1)\cap...\cap D(q_s)$]. Then $q$ is sectorial. If all $q_j$ are closed, so is $q$. If all the $q_j$ are closable so is $q$ and
$$\tilde{q}\subset \tilde{q_1}+...+\tilde{q_s}.$$
\end{thm}
The following theorem naturally associates  to strictly $m$-accretive quadratic forms a unique operator $T$.
\begin{thm} \label{thm:repstrcaccre} (\cite{RS1}, page 281) Let $q$ be a strictly $m$-accretive quadratic form with domain $Q(q)$. Then there is a unique operator $T$ on $\mathscr{H}$ such that:
\begin{itemize}
\item[a)] $T$ is closed.
\item[b)] $D(T)\subset Q(q)$ and if $\varphi,\psi \in D(T)$, then $q(\varphi,\psi)=(\varphi,T\psi)$. Further, $D(T)$ is a form core for $q$.
\item[c)] $D(T^*)\subset Q(q)$ and if  $\varphi,\psi \in D(T)$, then $q(\varphi,\psi)=(T^*\varphi,\psi)$.Further, $D(T^*)$ is a form core for $q$.
\end{itemize}
\end{thm}
From this theorem we can define.
\begin{defin}
A closed operator  $T$ is called {\bf strictly $m$-sectorial operator} if there exists $q$ strictly $m$-sectorial such that $q$ and $T$ satisfy properties a),b) c) of the above theorem.
\end{defin}
Now we can formulate the Ichinose lemma.
\begin{thm} (\cite{RS2},page 183) (Ichinose's lemma)\label{thm:Ichinoselemma} Let $\overline{S}_{\omega,\varphi,\theta}$ denote the sector $\{ z \vert \varphi-\theta \leq arg(z-\omega) \leq \varphi+\theta;\theta>\pi/2 \}$. Let $A$ and $B$ be strictly $m$-sectorial operators on Hilbert spaces $\mathscr{H}_1$ and $\mathscr{H}_2$ with sectors $\overline{S}_{\omega_1,\varphi,\theta_1}$ and $\overline{S}_{\omega_2,\varphi,\theta_2}$ (same $\varphi$!). Let $C$ denote the closure of $A\otimes I+I \otimes B$ on $D(A)\otimes D(B)$. Then $C$ is a strictly $m$-sectorial operator with sector $\overline{S}_{\omega_1+\omega_2,\varphi,min\{\theta_1,\theta_2\}}$ and $\sigma(C)=\sigma(A)+\sigma(B)$.
\end{thm}


\end{document}